\documentclass[11pt,a4paper]{article}

\usepackage{amsmath} % assumes amsmath package installed
\usepackage{amssymb}  % assumes amsmath package installed
\usepackage{graphics}
\usepackage{graphicx}
\usepackage{hyperref}
\hypersetup{
	colorlinks=true,
	linkcolor=black,
	filecolor=black,
	urlcolor=cyan,
	citecolor=black,
}
\usepackage{verbatim}
\usepackage{epsfig}
\usepackage{amsfonts}
\usepackage{mathrsfs}
\usepackage{bbm}
\usepackage{cprotect}
\usepackage{color}
\usepackage{makeidx}
\usepackage{color}
\usepackage{subfigure}
\usepackage{amsfonts}%
\usepackage{algorithm}
\usepackage{algorithmic}
\usepackage{amsthm}
\usepackage{nccmath}
\usepackage{bm}
\usepackage{multirow}
\usepackage{blkarray}
\usepackage{enumitem}
\usepackage{appendix}
\def\eps{\epsilon}

\DeclareMathOperator*{\argmin}{arg\,min}

\newcommand{\xs}{\hat{x}}
\newcommand{\ws}{\hat{w}}
\newcommand{\norm}[1]{\left\lVert#1\right\rVert}

\newtheorem{theorem}{Theorem}
\newtheorem{lemma}[theorem]{Lemma}

\makeindex

\title{Sparse data-driven quadrature rules via $\ell^p$-quasi-norm minimization} 
\author{Mattia Manucci\footnotemark[1] \and Jose Vicente Aguado \footnotemark[2] \and Domenico Borzacchiello \footnotemark[3]}

\begin{document}

\maketitle
\renewcommand{\thefootnote}{\fnsymbol{footnote}}

\footnotetext[1]{Gran Sasso Science Institute,
	             via Crispi 7,
	             L'Aquila, Italy. Email: {\tt mattia.manucci@gssi.it}}
             
\footnotetext[2]{$\acute{\text{E}}$cole Centrale de Nantes,
             	1 Rue de la No$\ddot{\text{e}}$,
             	Nantes, France. Email: {\tt Jose.Aguado-Lopez@ec-nantes.fr}}

\footnotetext[3]{$\acute{\text{E}}$cole Centrale de Nantes,
	             1 Rue de la No$\ddot{\text{e}}$,
                 Nantes, France. Email: {\tt domenico.borzacchiello@ec-nantes.fr}}

\renewcommand{\thefootnote}{\arabic{footnote}}

\begin{abstract}
In this paper we show the use of the focal underdetermined system solver to recover sparse empirical quadrature rules for parametrized integrals from existing data, consisting of the values of given parametric functions sampled on a discrete set of points. This algorithm, originally proposed for image and signal reconstruction, relies on an approximated $\ell^p$-quasi-norm minimization. The choice of $0<p<1$ fits the nature of the constraints to which quadrature rules are subject, thus providing a more natural formulation for sparse quadrature recovery compared to the one based on $\ell^1$-norm minimization.\\
We also extend an \textit{a priori} error estimate available for the $\ell^1$-norm formulation by considering the error resulting from data compression. Finally, we present two numerical examples to illustrate some practical applications. The first concerns the fundamental solution of the linear 1D Schrödinger equation, the second example deals with the hyper-reduction of a partial differential equation modelling a nonlinear diffusion process in the framework of the reduced basis method. For both the examples we compare our method with the one based on $\ell^1$-norm minimization and the one relaying on the use of the non-negative least square method. \texttt{Matlab} codes related to the numerical examples and the algorithms described are provided, see \cite{p19}.
\end{abstract}

{\bf Keywords: }Sparse quadrature rules, \texttt{FOCUSS} algorithm, linear programming, parametrized integrals, parametrized PDEs, hyper-reduction.


\pagestyle{myheadings}
\thispagestyle{plain}
\markboth{D.~Borzacchiello and M.~Manucci}{}

\section{Background and motivations}
The fundamental aim of the present work is to provide a reliable and efficient numerical approximation of integrals
\begin{equation}
I_k(\mu)=\int_{\Omega}f_k(x;\mu)\;\text{dx},\;\;\forall k\in\mathbb{K},\;\forall\mu\in\mathcal{D}.
\label{eq1.2}
\end{equation}
by means of \textit{sparse} empirical quadrature rules, for a given a family of functions 
\begin{equation}
f_k(x;\mu)\in L^{\infty}(\Omega),\;\Omega\subset\mathbb{R}^d,\;k\in \mathbb{K}=\{1,...,K\},\;\mu\in\mathcal{D}\subset\mathbb{R}^p.
\label{I.1}
\end{equation}

Parametrized integrals are found in many computer applications and their fast evaluation is often a key factor for performance optimization. For instance, such integrals are needed in transform methods for solving time-dependent ordinary differential equations, in which case $x$ is the frequency and $\mu$ is the time. Variational approximations schemes for partial differential equations (PDE) also involve integrals of this kind in which  $x$ is a space coordinate and $\mu$ can include physical parameters of the systems such as material or geometrical properties, external loads and source terms. In this case, computing inexpensive approximations of the integrals is of fundamental importance in the framework of reduced order methods relying on online-offline splitting of computational phases. \\ 

In general, the numerical approximation of integrals is performed by means of quadrature rules such as
\begin{equation}
I_k(\mu)=\int_{\Omega}f_k(x;\mu)\;\text{dx}\approx\sum_{i=1}^{\mathcal{N}}w_if_k(x_i;\mu).
\end{equation} 
The accuracy of the approximation depends on the sampling strategy adopted for the selection of quadrature nodes $x_i$. Even when considering complex domains of integration, it is usually possible to achieve high accuracy using composite quadrature rules. These are largely used in piece-wise approximations for the numerical solution of PDEs. We refer to these kind of quadrature rules as \textit{full order quadrature rules} and denote them by the pair of nodes and weights $\{x,w\}$. If high accuracy is desired, $\mathcal{N}$ is expected to grow large, resulting in an elevated cost for numerical integration. 

The starting point to build efficient quadrature rules is the consideration that, if $f_k(\mu)$ lies in a lower dimensional manifold, it is possible to perform numerical integration using a considerably smaller number of nodes, $\mathcal{K}\ll\mathcal{N}$, without significantly lowering the accuracy prescribed for the full order rule. This is often the case of the applications mentioned before, in which parametric integrals have to be evaluated multiple times, that is $\forall k\in\mathbb{K}$ and $\forall \mu\in\mathcal{D}$.  We refer to these new rules as \textit{sparse} or \textit{reduced} quadrature rules and use the notation $\{\xs,\ws\}$. 

Sparse quadrature rules can be calibrated by solving a sparse linear regression problem on the paired data consisting of function values $f_k(x_i,\mu_m)$ sampled at nodes $x_i$, $i=1,...,\mathcal{N}$ for some $\mu_m$, $m\in\mathbb{M}=\{1,...,M\}$ and the corresponding full order approximation of the integrals $I_k(\mu_m)$.\\

The use of sparsity promoting techniques implies that the nodes of the new rule $\xs_i$, $i=1,...,\mathcal{K}$ are opportunely selected among those of the full order quadrature rule. In the related literature this kind of rules are also called \textit{empirical quadrature rule} since their derivation is form a computed dataset. This approach has already been adopted in a previous work of the authors and discussed in more details in \cite{p12}.\\

The core differences between the methodologies developed to recover empirical rules concern either the formulation of the associated the regression problem or the numerical technique used to solve it. One possible approach is to find an approximation by interpolation for $f_k$ first and then to perform the integration as a linear combination of function values at interpolation points through the integrals of the interpolation functions, which can be regarded as quadrature weights and can be computed during the off-line phase. This techniques are  called ``interpolate-then-integrate". In this case, the problem reduces to finding a set of optimal interpolation points for the family of functions $f_k$. If the parametric functions belong to a low-dimensional subspace, these points can be found through the Empirical Interpolation Method \cite{p6}. The discrete counterpart of this method is named DEIM \cite{p15} and can operate directly with discrete functions. The main drawback of this approach is that the quadrature weights are not guaranteed to be positive which can lead to stability issues of the numerical schemes for the discretization of PDEs.\\

The procedures aiming at the straightforward approximation of the integral are also referred to as \textit{empirical cubature methods}, as in \cite{p7} and \cite{p8}. These methods works in the framework of residual minimization in $\ell^2$-norm, therefore sparsity has to be explicitly enforced, for example through a heuristic sequential point selection process \cite{p7} or through an approximate $\ell^0$-pseudo-norm penalization term, \cite{p8}. In both cases, the main bottleneck is the solution of non-negative least-square problems arising from the non negativity constraint that is imposed to the integration weights. This can be potentially expensive when applied to large datasets. We also mention the work in \cite{p18}, where it is shown how to partially avoid the solution of non-negative least-square problem. Another way to obtain sparse vectors of non-negative weights was proposed by Ryu and Boyd in \cite{p9}. It consists in replacing the $\ell^0$-pseudo-norm minimization with the $\ell^1$-norm. This norm naturally yields quadrature rules that are sparse and furthermore the offline problem can be cast as a linear program (LP) and efficiently treated by the SIMPLEX algorithm \cite{p10}. The linear programming methodology to recover sparsity has been employed by Patera and Yano \cite{p1,p4} producing a method that, for a desired accuracy $\delta$ and a large enough training dataset, is able to provide a reduced rule with an integration error that is within the prescribed tolerance $\delta$. The use of $\ell^1$-norm minimization of the weight vector to achieve sparsity is somehow counterintuitive, since for non negative weights the $\ell^1$-norm is constrained to be equal to the measure of the domain of integration.  Although this approach sets out what seem to be competing objectives, very reasonable results can be obtained by simply relaxing the imposition of the constraint (i.e. allowing some tolerance in its fulfillment). According to our numerical experiments, only when a strict tolerance is demanded the method fails to achieve the target accuracy.\\
The alternative solution proposed in this work relies on the minimization of the $\ell^p,\;0<p<1$, quasi-norm\footnote{Is said to be a quasi-norm since it violates the triangular inequality for $p\in(0,1)$.} minimization as way to enforce sparsity. In this way the constraint on the sum of the weights can be enforced exactly. In the same line of the work presented in \cite{p1}, we are able to prove that our method recovers sparse rules which integrate functions with an error proportional to a prescribed accuracy $\epsilon$ even when $\epsilon$ is very small.

\section{Problem formulation} \label{sec:pformulation}

The numerical approximation of integrals in (\ref{eq1.2}) with the use of a full quadrature rule is defined by
\begin{equation}
I_k(\mu)=\int_{\Omega}f_k(x;\mu)\;\text{dx}, \;\;\;\;\;I_k^{\text{full}}(\mu)=\sum_{i=1}^{\mathcal{N}}w_if_k(x_i;\mu).
\end{equation}
The rule is required to yield an accurate approximation of the integral under a prescribed tolerance $\varepsilon/2$ for the integration error 
\begin{equation}
|I_k(\mu)-I_k^{\text{full}}(\mu)|\le\varepsilon/2,\;\forall k=1,...,K,\;\forall\mu\in\mathcal{D};
\end{equation}
with additional constraints regarding the exact integration of the constant function and positivity of the quadrature weights, i.e.
\begin{align}
&\sum_{i}^{\mathcal{N}}w_i=|\Omega|,\label{eq1.4}\\
&\;w_{i}\ge0,\;\forall i;\label{eq1.10}
\end{align} 
where $|\Omega|$ stays for the the measure of the domain. As already mentioned, for strictly non negative quadrature weights, (\ref{eq1.4}) expresses a constraint on the $\ell^1$-norm of the weight vector. We aim to find a sparse rule which does not deteriorate the quality of the approximation and satisfies constraints (\ref{eq1.4}) and (\ref{eq1.10}), i.e.
\begin{align}
&I^{\text{sparse}}_k(\mu)= \sum_{i=1}^{\mathcal{K}}\ws_if_k(\xs_i;\mu),\;\mathcal{K}\ll\mathcal{N}\\
&|I^{\text{full}}_k(\mu)-I_k^{\text{sparse}}(\mu)|\le\varepsilon/2,\;\forall k=1,...,K,\;\forall\mu\in\mathcal{D};\label{eq1.8}\\
&\sum_{i}^{\mathcal{N}}\ws_i=|\Omega|,\;\label{eq1.8bis}\\
&\ws_{i}\ge0,\;\forall i;\label{eq1.7}
\end{align}
so that the final integration error given is
\begin{equation}
|I_k(\mu)-I_k^{\text{sparse}}(\mu)|\le|I_k(\mu)-I_k^{\text{full}}(\mu)|+|I^{\text{full}}_k(\mu)-I_k^{\text{sparse}}(\mu)|\le\varepsilon.
\end{equation}
To find such a sparse quadrature rule we consider a training dataset $\Xi_{\mathbb{M}}^{\text{train}}=\{\mu_m^{\text{train}}\in\mathcal{D}\}_{m\in\mathbb{M}}$, $|\Xi_{\mathbb{M}}^{\text{train}}|=N^{\text{train}}$, and we approximate the integrals (\ref{eq1.2}) for each $\mu\in\Xi_{\mathbb{M}}^{\text{train}}$. Considering also the constraint (\ref{eq1.4}), the number of integrals to compute is $K\cdot N^{\text{train}}+1$. The problem can be cast in a matrix-vector form as
\begin{equation}
Aw=b;
\label{eq1.3}
\end{equation} 
where $A\in\mathbb{R}^{(K\cdot N^{\text{train}}+1)\times\mathcal{N}}$, $w\in\mathbb{R}^{\mathcal{N}}$ and $b\in\mathbb{R}^{K\cdot N^{\text{train}}+1}$. Column $j$ of $A$ corresponds to the evaluations $f_k(\;\cdot\;;\mu_i)\;\forall \;k,\;i,$ in the quadrature point $x_j$ while the last row is a vector of one corresponding to the constant function, vector $w$ contains the weights of the full quadrature rule while each entry of $b$ is equal to $I_k^{\text{full}}(\mu_{j})$ and the last entry is equal to $|\Omega|$.\\
Assuming $(K\cdot N^{\text{train}}+1)<\mathcal{N}$ and taking the quadrature weights $w$ as unknowns, we can view (\ref{eq1.3}) as an underdetermined system, i.e. a  system with possibly infinite many solutions. Among them we look for those solutions with least non-zero entries, the so called \textit{sparse solutions}, so the problem can read as
\begin{equation}
\min_{y\in\mathbb{R}^{\mathcal{N}}}\|y\|_{\ell^0},
\label{eq1.5}
\end{equation}
subject to
\begin{align}
&Ay=b,\;\label{eq1.6}\\
&y_i\ge 0,\;\forall i.\label{eq1.6bis}
\end{align}
Note that (\ref{eq1.3}) is not formally equivalent to the restriction expressed by  (\ref{eq1.8}) on the train sample. Indeed, it is a more restrictive condition, since it expresses a set of equations rather than inequalities. The  matrix-vector form is introduced here to formulate the set of constraints as an underdetermined system. However, it will be clear from the next sections that the convergence of the iterative solution of the system of equations  (\ref{eq1.5})-(\ref{eq1.6})-(\ref{eq1.6bis})  within the given tolerance implies that 
\begin{equation}\label{eq:ineq}
	\|Ay-b\|\le \varepsilon/2.
\end{equation}
Therefore, if (\ref{eq:ineq}) holds then (\ref{eq1.8}) holds $\forall \mu \in \Xi^{\text{train}}_\mathbb{M}$.

\section{From $\ell^0$-pseudo-norm to $\ell^p$-quasi-norm minimization problem}
The $\ell^0$-pseudo-norm minimization (\ref{eq1.5}) needs an appropriate reformulation to reduce its complexity as the problem is NP hard, as stated in \cite{p17}. In \cite{p1,p4,p9} the $\ell^1$-norm is chosen as surrogate of the $\ell^0$-pseudo-norm since it naturally provides sparse solution when combined with the constraints (\ref{eq1.6}). Moreover, the obtained formulation can then be cast as a linear programming problem for which efficient algorithms are available.\\
We propose a new approach where the $\ell^0$-pseudo-norm is replaced by the $\ell^p$-quasi-norm for $0<p<1$. The use of the $\ell^p$-quasi-norm minimization can be motivated by the following rationale: the constraints in (\ref{eq1.7}) are equivalent to  
\begin{equation}
\|y\|_{\ell^1}=|\Omega|;
\label{eq1.9}
\end{equation} 
which is the quantity that the problem seeks to minimize. This creates a contrast between the objective functions as promoting sparsity conflicts with one of the accuracy constraints. As  shown in \cite{p2,p3}, if $\ell^1$-norm minimization succeeds in the optimization problem (\ref{eq1.5}), then also $\ell^p$-quasi-norm for any $0<p<1$ is able to recover a sparse solution of (\ref{eq1.3}). Therefore, the $\ell^p$-quasi-norm seems a more suitable approach for sparse constrained minimization in this specific case.\\
We must mention that in \cite{p1} constraints (\ref{eq1.6}) are required to be  satisfied only up to an accuracy $\epsilon_1$. This relaxation comes naturally since by (\ref{eq1.8}) the new sparse rule is considered as an approximation of the full order rule. Since (\ref{eq1.9}) is not exactly satisfied,  the conflict with the objective function is somehow resolved, though, for high accuracy requirements, it can induce a deterioration in the approximation.\\
In the light of these considerations, we can formulate the problem in the following way:
\begin{align}
&\min_{y\in\mathbb{R}^{\mathcal{N}}}\|y\|_{\ell^p},\label{eq2.6}\\
&\|Ay-b\|_2\le\epsilon_1,\; \label{eq2.5}\\
&y_i\ge 0,\;\forall i.\label{eq2.8}
\end{align}
where the $\ell^p$-quasi-norm is formally defined by
\begin{equation}
	\norm{y}_{\ell^p}=\left(\sum_{i}|y_i|^p\right)^{\frac{1}{p}}.
\end{equation}
\subsection{The FOCUSS algorithm} \label{subsec:FOCUSS alorithm}
We solve the $\ell^p$-quasi-norm optimization problem by using the \textit{Focal Underdetermined System Solver algorithm} or \verb|FOCUSS| \cite{p5}. This method is an iterative fixed-point algorithm that achieves the $\ell^p$-quasi-norm minimum by a sequence of weighted $\ell^2$-norm optimization problems, each having a unique solution. Given the solution $y_k$ at step $k$ the solution at step $k+1$ is computed as
\begin{align}
&y_{k+1}=W_k(AW_k)^{\dagger}b,\nonumber\\
&W_{k}=\text{diag}(y^q_k);
\end{align}
where $(AW_k)^{\dagger}$ is the Monroe-Penrose pseudoinverse. After some iterations of the algorithm we have that $y_{k-1}\approx y_k$ and consequentially the objective minimized at each step becomes
	\begin{equation*}
		||W_{k-1}y_k||^2_{2}=\sum_{i=1,(y_i)_{k-1}\ne 0}^{\mathcal{N}}\left(\frac{(y_i)_k}{(y_i)^{q}_{k-1}}\right)^2\approx\sum_{i=1}^{\mathcal{N}}(y_i)_k^{2-2q}
	\end{equation*}  
showing that the successive minimization of the weighted $\ell^2$-norm leads to the minimization of the norm $\ell^{2-2q}$. The further requirement that $0<2-2q<1$ leads to  the constraint $0.5<q<1$.\\
At every step $k$, the solution $y_k$ satisfies exactly\footnote{Up to machine precision.}, i.e $\epsilon_1=0$, the constraints (\ref{eq2.5}). Since the full rule has a known integration error, we accept the same order of error from the empirical rule. Moreover, we expect that this relaxation of constraints (\ref{eq1.6}) into (\ref{eq2.5}) provides a solution with a smaller number of nonzero entries.\\
This idea is implemented through regularization, which consists in changing the objective function (\ref{eq2.6}) into 
\begin{equation}
\argmin_{y\in\mathbb{R}^{\mathcal{N}}}\left(\|b-Ay\|^2_2+\lambda\|y\|^p_p\right),
\label{eq2.7}
\end{equation}
where $\lambda\in\mathbb{R}^{+}$ is a parameter that can be adjusted so as to prioritize sparsity over accuracy (for $\lambda\gg1$) or vice-versa. Note that for $p=2$ this is the well known Tikhonov regularization.\\
The corresponding iterative scheme for the minimization of the regularized problem, as presented in \cite{p5}, is: 
\begin{equation}
y_{k+1}=W_kW_kA^{T}(AW_kW_kA^{T}+\lambda I)^{-1}b.
\label{eq1.11}
\end{equation}
For the remainder of this paper, we assume that at each iteration the solution $y_k$ satisfies (\ref{eq2.5}). In the next subsection we discuss how to modify the algorithm so that this assumption is always verified and in particular we will show how to get a solution for which the $\ell^2$-norm of the residual is exactly equal to $\epsilon_1$.\\
\subsection{The FOCUSS algorithm for sparse quadrature recovery}
In order to adapt the \verb|FOCUSS| algorithm to our problem, the original version needs to ensure the non negativity of the integration weights. As already discussed in \cite{p5}, a solution consists in introducing a relaxation step that is performed after each iteration.\\
In practice, this consists in finding $\alpha_k\in(0,1)$ such that
\begin{align}
&y_i^{\text{new}}\ge 0,\; \forall i=1,...,\mathcal{N},\nonumber\\
&y_k^{\text{new}}=\alpha_ky_k+(1-\alpha_k)y_{k-1}.\label{3.1}
\end{align}
We start from an initial guess with all positive entries, for instance the full order quadrature rule, and then at every iteration $k$ we check that each entry is still positive. If not, we determine the new solution of step $k$ according to the condition expressed in (\ref{3.1}).\\
At each iteration $k+1$ of \verb|FOCUSS|, the pseudoinverse of the matrix $AW_k$ needs to be computed. Here,  $A$ is a matrix of $\mathcal{R}$ rows, with $\mathcal{R}\ll\mathcal{N}$, and $W_k$ is the diagonal matrix of weights. This step requires the solution of the linear system $[AW_k(AW_k)^{T}]c=b$, for which the size of the coefficient matrix is $\mathcal{R}\times \mathcal{R}$. Thus the computational cost is given by the matrix multiplication plus the solution of the system, that is, $2\mathcal{N}\mathcal{R}^2+\frac{\mathcal{R}^3}{6}\approx \rm{O}(\mathcal{N}\mathcal{R}^2)$. A way to compute the solution associated with the pseudoinverse is given by the truncated singular values decomposition. If $[U_k,S_k,V_k]=\text{svd}(AW_k)$ then it can be shown that
\begin{align*}
y_{k+1}=W_{k}[AW_k]^{\dagger}b&=W_{k}V_kS_kU_k^{T}[U_kS_kV_k^{T}V_kS_kU_k^{T}]^{-1}b\\
&=W_{k}V_kS_kU_k^{T}[(U^{T}_k)^{-1}(U_kS_k^2)^{-1}]b\\
&=W_{k}V_kS_k^{-1}U_k^{T}b=\sum_{n=1}^{\mathcal{R}}\frac{<u_{n}^T,b>_n^q}{\sigma_n}v_n*y^q_{k},
\end{align*}
where $<\cdot,\cdot>$ denotes the Euclidean inner product and * stands for the component wise product of two vectors of same size. Note that the computational cost of SVD by direct algorithms is of the order $\rm{O}(\mathcal{N}\mathcal{R}^2)$.\\
We recall that our aim is to solve the Tikhonov regularization problem. Therefore our solution at step $k+1$ is given by (\ref{eq1.11}), for which we know the singular values decomposition of $AW_k$ and we obtain
\begin{align*}
y^r_{k+1}=&W_kW_kA^{T}(AW_kW_kA^{T}+\lambda I)^{-1}b\\
=&W_kV_kS_kU_k^{T}[U_kS^2_kU_k^{T}+\lambda I]^{-1}b \\
=&W_kV_kS_kU_k^{T}[U_k(S^2_k+\lambda I)U_k^{T}]^{-1}b\\
=&W_kV_kS_k[S^2_k+\lambda I]^{-1}U_k^{T}b\\
=&\sum_{n=1}^{\mathcal{R}}\frac{\sigma_n}{\sigma^2_n+\lambda}<u_{n}^T, b>v_n*y^q_{k}.
\end{align*}
The key point is the control of the residual at each iteration $k+1$, that is:
\begin{equation*}
\|b-Ay^r_{k+1}\|_2\le\epsilon_1,\;\forall\;k.
\end{equation*}
Once one has the SVD decomposition of $AW_k$ the evaluation of the norm of the residual is immediate, indeed
\begin{align}
\|b-Ay^r_{k+1}\|_2=&\|A(y_{k+1}-y^r_{k+1})\|_2\nonumber\\
=&\|AW_kV_k[S_k^{-1}-S_k(S^2_k+\lambda I)^{-1}]U_k^Tb\|_2\nonumber\\
=&\|U_kS_k[S_k^{-1}-S_k(S^2_k+\lambda I)^{-1}]U_k^Tb\|_2\nonumber\\
=&\Bigg\|\sum_{n=1}^{\mathcal{R}}\frac{\lambda}{\sigma_n^2+\lambda}<u_n^T,b>u_n\Bigg\|\nonumber\\
=&\sqrt{\sum_{n=1}^{\mathcal{R}}\left(\frac{\lambda}{\sigma_n^2+\lambda}<u_n^T,b>\right)^2}.
\label{eq3.6}
\end{align}
The expression found in (\ref{eq3.6}) is quite inexpensive and can be easily evaluated for different values of $\lambda$ once the truncated SVD decomposition of $AW_k$ is available. Therefore, a possible strategy can be to compute the truncated SVD factorization of $AW_k$ and then find $\lambda$ such that $\|b-Ay_{k+1}\|_2=\epsilon_1$ at each iteration $k+1$. In this way it is possible to reduce the number of non-zero weights for a prescribed accuracy without increasing the order of the computational cost.\\
The process is iterated until the following three conditions are simultaneously satisfied:
	\begin{enumerate}
		\item $\|y_k-y_{k-1}\|\le tol$, with $tol$ a prescribed threshold;
		\item $\|b-Ay_k\|\le\epsilon_1$;
		\item $\mathcal{K}\le\mathcal{R}$.
	\end{enumerate}
In terms of algorithm iterations, since condition $1$ is typically verified later than conditions $2$ and $3$, we provide an additional stagnation criterion on $\mathcal{K}$. If after $Q\in\mathbb{N}$ iterations $\mathcal{K}$ is unchanged, the algorithm is stopped.
\subsection{The truncated Singular Value Decomposition for data compression} \label{sec:POD}
In practical applications, it is often likely that $K\times N^{\text{train}}+1$ is a large number, e.g. we have a large training dataset and/or a large family of functions, producing a matrix $A$ with rank smaller than $K\times N^{\text{train}}+1$. This is the case of a dataset that contains redundant information which generates additional computational cost and possibly  translates in an ill conditioned problem.\\
To avoid this situation, we perform a truncated Singular Value Decomposition (SVD) of $A^{T}$: we extract the most significant modes, up to an established tolerance, and then we run the \verb|FOCUSS| algorithm.\\
Having that $A^{T}=USV^{T}$ and considering the spectral energy up to mode $\mathcal{R}$, the system of constraints becomes
\begin{equation}
\underbrace{S(1:\mathcal{R},1:\mathcal{R})U(\;:\;,\;1:\mathcal{R})^{T}}_{\tilde{A}}y=\underbrace{V(\;:\;,\;1:\mathcal{R})^{T}b}_{\tilde{b}},
\label{eq14}
\end{equation}
where the notation $S(1:\mathcal{R},1:\mathcal{R})$ indicates that we extract from matrix $S$ the first $\mathcal{R}$ rows and columns while $U(\;:\;,\;1:\mathcal{R})$ means that we select from $U$ all the rows and the first $\mathcal{R}$ columns.\\
Naturally, with this approach, we are introducing the error due to the neglected modes. However it is well known how to quantify this error and in particular it is known how to select $\mathcal{R}$ in such a way that the error introduced is smaller than a prescribed accuracy $\epsilon_2$.\\
Given a training dataset $\Xi^{\text{train}}_{\mathbb{M}}$ we define the set of snapshots $\{\phi_{k,m}\equiv f_k(\;\cdot\;;\mu_m)\}_{k\in\mathbb{K},m\in\mathbb{M}}$. Note that for a simplification of the notation in this paragraph we are not considering constrain (\ref{eq1.9}). From the dataset we perform the truncated SVD and we extract the $r$ most representative modes. The error $E^{\mathcal{R}}_{SVD}$ produced by the SVD base $\zeta_1,...,\zeta_\mathcal{R}$ of dimension $\mathcal{R}$ in the approximation of the entire set of snapshots $\{\phi_{k,m}\}_{k\in\mathbb{K},m\in\mathbb{M}}$, defined as
\begin{align}
E^{\mathcal{R}}_{SVD}=&E(\zeta_1,...,\zeta_\mathcal{R})=\sum_{m=1}^{N^{\text{train}}}\sum_{k=1}^{K}||\phi_{k,m}-\Pi_{\zeta_\mathcal{R}}\left[\phi_{k,m}\right]||_2^2,\\
&\Pi_{\zeta_\mathcal{R}}\left[\phi\right]=\sum_{n=1}^{\mathcal{R}}\pi_{\zeta_\mathcal{R}}^{n}\left[\phi\right]\zeta_n,\;\text{with}\;\pi_{\zeta_\mathcal{R}}^{n}\left[\phi\right]=<\phi,\zeta_n>;\label{2.1.1}
\end{align} 
is equal to the sum of the square of the singular values
\begin{equation}
E(\zeta_1,...,\zeta_\mathcal{R})=\sum_{i=\mathcal{R}+1}^{K\cdot N^{\text{train}}}\sigma_i^2,
\label{eq9}
\end{equation}
related to the $K\cdot N^{\text{train}}-\mathcal{R}$ modes that were not selected for the basis. Therefore it is enough to choose $\mathcal{R}$ as the smallest value of $\tilde{\mathcal{R}}$ such that
\begin{equation}
I(\tilde{\mathcal{R}})=\frac{\sum_{i=1}^{\tilde{\mathcal{R}}}\sigma_i^2}{\sum_{i=1}^{K\cdot N^{\text{train}}}\sigma_i^2}\le1-\epsilon_2.
\end{equation} 
Once the basis $\{\zeta_{i}\}_{i=1,...,\mathcal{R}}$ is determined the formulation in (\ref{eq2.6}), (\ref{eq2.5}), (\ref{eq2.8}) can be modified as follow: find a vector $y\in\mathbb{R}^{\mathcal{N}}$ such that
\begin{align}
\text{minimize}& \left(\sum_{i=1}^{\mathcal{N}}y_i^{\text{p}}\right)^{1/p}\label{2.4.2}\\
\text{subject to }&\nonumber\\
&\Big\|\tilde{A}y-\tilde{b}\Big\|\le\epsilon_1,\label{2.4.1}\\
\text{and}\;\;&y_i\ge0,\;1\le i\le\mathcal{N}.\label{2.4.3}
\end{align}
The solution of this problem is found using \verb|FOCUSS| with Tikhonov regularization. We then identify the indices associated with the non-zero values of $y$ as $i_k$, $1\le k\le \mathcal{K}$, and set $\xs_k=x_{i_k}$, $\ws_k=y_{{i_k}}$, $1\le k\le \mathcal{K}$. This approach also provides an upper bound for the number of non-zero entries of $y$, that approaches $\mathcal{R}$ as $\epsilon_1\rightarrow0$.
\section{Error Analysis}  \label{sec: Error Analysis}
In this section we build on the work presented in \cite{p1} and add an additional term to the error analysis due to the compression of the data through truncated SVD.\\
For $\mathcal{K}$ being the number of non-zero entries of $y$, we have that $\mathcal{K}\le \mathcal{R}\ll \min\{K\cdot N^{\text{train}}{+1},\mathcal{N}\}$. For this to hold, it must be possible to construct a low rank approximation of the manifold 
\begin{equation}
\mathcal{M}=\{f_k(\mu):\;k\in\mathbb{K},\;\mu\in\mathcal{D}\},
\end{equation}
by applying a dimensionality reduction technique to the training dataset.\\
We now provide a general results quantifying the error associated to the procedure proposed in this paper.\\
\begin{lemma}
	For any $\mu\in\mathcal{D}$ it holds

		\begin{align}
		&\max_{k\in\mathbb{K}}|I_k^{\text{full}}(\mu)-I_m^{\text{sparse}}(\mu)|\nonumber\\
		&\le\max_{k\in\mathbb{K}}\bigg(\inf_{\alpha\in\mathbb{R}^{N^{\text{train}}}}\bigg(\left(\|w\|_{2}+\|\ws\|_{2}\right)\sqrt{\sum_{m=1}^{N^{\text{train}}}\alpha_m^2}\sqrt{\sum_{i=\mathcal{R}+1}^{K\cdot N^{\text{train}}+1}\sigma_i^2}\label{eqLemma}\\
		&+\epsilon_1\max_{m\in\mathbb{M}}\left(\sum_{n=1}^{\mathcal{R}}|\pi_{\zeta_R}^{n}\left[\phi_{k,m}\right]|\right)\sum_{m=1}^{N^{\text{train}}}|\alpha_m|+2|\Omega|\Big\|f_k(\;\cdot\;;\mu)-\sum_{m=1}^{N^{\text{train}}}\alpha_m\phi_{k,m}\Big\|_{L^{\infty}(\Omega)}\bigg)\bigg)\nonumber.
		\end{align}
	 
	\label{lemma1}
\end{lemma}
\begin{proof}
	First, we fix $k\in\mathbb{K}$. Then, for any $\alpha\in\mathbb{R}^{N^{\text{train}}}$ we obtain
	\begin{align}
	|I_k^{\text{full}}&(\mu)-I_k^{\text{sparse}}(\mu)|=\Bigg|\sum_{i=1}^{\mathcal{N}}w_{i}f_k(x_{i};\mu)-\sum_{i=1}^{\mathcal{K}}\ws_{i}f_k(\xs_{i};\mu)\Bigg|\nonumber\\
	\le&\Bigg|\sum_{i=1}^{\mathcal{N}}w_{i}\sum_{m=1}^{N^{\text{train}}}\alpha_m\phi_{k,m}(x_{i})-\sum_{i=1}^{\mathcal{K}}\ws_{i}\sum_{m=1}^{N^{\text{train}}}\alpha_m\phi_{k,m}(\xs_{i})\Bigg|\nonumber\\
	+&\Bigg|\sum_{i=1}^{\mathcal{N}}w_i\left(f_k(x_{i};\mu)-\sum_{m=1}^{N^{\text{train}}}\alpha_m\phi_{k,m}(x_{i})\right)\Bigg|\label{eq12}\\
	+&\Bigg|\sum_{i=1}^{\mathcal{K}}\ws_i\left(f_k(\xs_{i};\mu)-\sum_{m=1}^{N^{\text{train}}}\alpha_m\phi_{k,m}(\xs_{i})\right)\Bigg|,\nonumber
	\end{align}
	where the second term follows by H\"older's inequality ($p=1,\;q=\infty$) and (\ref{eq1.4})-(\ref{eq1.10}), and the third by using the H\"older's inequality and recalling that $\sum_{k=1}^{\mathcal{K}}\ws_k\le\sum_{i=1}^{\mathcal{N}}w_{i}=|\Omega|$. On the other hand, considering the first term we have that
	%For the second term in (\ref{eq12}) we first apply H\"older's inequality ($p=1,\;q=\infty$) and then we use (\ref{eq1.4})-(\ref{eq1.10}). For the third term we again employ H\"older's inequality and the fact that $\sum_{k=1}^{\mathcal{K}}w_k^{\text{sparse}}\le\sum_{i=1}^{\mathcal{N}}w_{i}^{\text{truth}}=|\Omega|$. We now focus on the first term,
	\begin{align}
	&\Bigg|\sum_{i=1}^{\mathcal{N}}w_{i}\sum_{m=1}^{N^{\text{train}}}\alpha_m\phi_{k,m}(x_{i})-\sum_{i=1}^{\mathcal{K}}\ws_{i}\sum_{m=1}^{N^{\text{train}}}\alpha_m\phi_{k,m}(\xs_{i})\Bigg|\nonumber\\
	=&\Bigg|\sum_{m=1}^{N^{\text{train}}}\alpha_m\left(\sum_{i=1}^{\mathcal{N}}w_{i}\phi_{k,m}(x_{i})-\sum_{i=1}^{\mathcal{K}}\ws_{i}\phi_{k,m}(\xs_{i})\right)\Bigg|\nonumber\\
	\le&\underbrace{\Bigg|\sum_{m=1}^{N^{\text{train}}}\alpha_m\left(\sum_{i=1}^{\mathcal{N}}w_{i}\phi_{k,m}(x_{i})-\sum_{i=1}^{\mathcal{N}}w_{i}\Pi_{\zeta_\mathcal{R}}\left[\phi_{k,m}\right](x_{i})\right)\Bigg|}_{\bf A}\nonumber\\
	+&\underbrace{\Bigg|\sum_{m=1}^{N^{\text{train}}}\alpha_m\left(\sum_{i=1}^{\mathcal{N}}w_{i}\Pi_{\zeta_\mathcal{R}}\left[\phi_{k,m}\right](x_{i}))-\sum_{i=1}^{\mathcal{K}}\ws_{i}\Pi_{\zeta_\mathcal{R}}\left[\phi_{k,m}\right](\xs_{i})\right)\Bigg|}_{\bf B}\label{eq10}\\
	+&\underbrace{\Bigg|\sum_{m=1}^{N^{\text{train}}}\alpha_m\left(\sum_{i=1}^{\mathcal{K}}\ws_{i}\Pi_{\zeta_\mathcal{R}}\left[\phi_{k,m}\right](\xs_{i})-\sum_{i=1}^{\mathcal{K}}\ws_{i}\phi_{k,m}(\xs_{i})\right)\Bigg|}_{\bf C}\nonumber.
	\end{align}
	In the following, we estimate the terms {\bf A}, {\bf B} and {\bf C} in \eqref{eq10} separately.
	First, by H\"older's inequality ($p,q=2$) for vectors in $\mathbb{R}^{\mathcal{N}}$ we have that
	\begin{align*}
	{\bf A} \le&\sum_{m=1}^{N^{\text{train}}}|\alpha_m|\Bigg|\sum_{i=1}^{\mathcal{N}}w_i\{\phi_{k,m}(x_{i})-\Pi_{\zeta_\mathcal{R}}\left[\phi_{k,m}\right](x_{i})\}\Bigg|
	\\
	\le&\sum_{m=1}^{N^{\text{train}}}|\alpha_m|\|w\|_{2,\mathbb{R}^{\mathcal{N}}}\|\phi_{k,m}-\Pi_{\zeta_\mathcal{R}}\left[\phi_{k,m}\right]\|_{2,\mathbb{R}^{\mathcal{N}}}\
	\end{align*}
	and, again by H\"older ($p,q=2$) for vectors in $\mathbb{R}^{N^{\text{train}}}$ we obtain
	\begin{align*}
	& \sum_{m=1}^{N^{\text{train}}}|\alpha_m|\|w\|_{2,\mathbb{R}^{\mathcal{N}}}\|\phi_{k,m}-\Pi_{\zeta_\mathcal{R}}\left[\phi_{k,m}\right]\|_{2,\mathbb{R}^{\mathcal{N}}}
	\\
	\le&\ \|w\|_{2,\mathbb{R}^{\mathcal{N}}}\sqrt{\sum_{m=1}^{N^{\text{train}}}\alpha_m^2}\sqrt{\sum_{m=1}^{N^{\text{train}}}\|\phi_{k,m}-\Pi_{\zeta_\mathcal{R}}\left[\phi_{k,m}\right]\|_{2,\mathbb{R}^{\mathcal{N}}}^2},
	\end{align*}
	where the notation $\|\cdot\|_{2,\mathbb{R}^{\mathcal{N}}}$ stands for the Euclidean norm on $\mathbb{R}^{\mathcal{N}}$. Then, by using \eqref{eq9} we get 
	\begin{align*}
	& \|w\|_{2,\mathbb{R}^{\mathcal{N}}}\sqrt{\sum_{m=1}^{N^{\text{train}}}\alpha_m^2}\sqrt{\sum_{m=1}^{N^{\text{train}}}\|\phi_{k,m}-\Pi_{\zeta_\mathcal{R}}\left[\phi_{k,m}\right]\|_{2,\mathbb{R}^{\mathcal{N}}}^2}
	\\
	&\le\|w\|_{2,\mathbb{R}^{\mathcal{N}}}\sqrt{\sum_{m=1}^{N^{\text{train}}}\alpha_m^2}\sqrt{\sum_{i=\mathcal{R}+1}^{K\cdot N^{\text{train}}+1}\sigma_i^2},
	\end{align*}
	which implies that 
	\begin{equation}\label{eq12bis}
	{\bf A} \leq \|w\|_{2,\mathbb{R}^{\mathcal{N}}}\sqrt{\sum_{m=1}^{N^{\text{train}}}\alpha_m^2}\sqrt{\sum_{i=\mathcal{R}+1}^{K\cdot N^{\text{train}}+1}\sigma_i^2}.	
	\end{equation}
	We now go to the analysis of term {\bf B}. By \eqref{2.1.1} we obtain
	\begin{equation*}
	{\bf B} = \Bigg|\sum_{m=1}^{N^{\text{train}}}\alpha_m\left(\sum_{i=1}^{\mathcal{N}}w_{i}\sum_{n=1}^{\mathcal{R}}\pi_{\zeta_\mathcal{R}}^{n}\left[\phi_{k,m}\right]\zeta_n(x_{i}))-\sum_{i=1}^{\mathcal{K}}\ws_{i}\sum_{n=1}^{\mathcal{R}}\pi_{\zeta_\mathcal{R}}^{n}\left[\phi_{k,m}\right]\zeta_n(\xs_{i})\right)\Bigg|	
	\end{equation*}
	and, by definition of $\ws_i$ and $\xs_i$ we deduce that
	\begin{align*}
	& 	\Bigg|\sum_{m=1}^{N^{\text{train}}}\alpha_m\left(\sum_{i=1}^{\mathcal{N}}w_{i}\sum_{n=1}^{\mathcal{R}}\pi_{\zeta_\mathcal{R}}^{n}\left[\phi_{k,m}\right]\zeta_n(x_{i}))-\sum_{i=1}^{\mathcal{K}}\ws_{i}\sum_{n=1}^{\mathcal{R}}\pi_{\zeta_\mathcal{R}}^{n}\left[\phi_{k,m}\right]\zeta_n(\xs_{i})\right)\Bigg|
	\\
	=&\Bigg|\sum_{m=1}^{N^{\text{train}}}\alpha_m\left(\sum_{n=1}^{\mathcal{R}}\pi_{\zeta_\mathcal{R}}^{n}\left[\phi_{k,m}\right]\left(\sum_{i=1}^{\mathcal{N}}w_{i}\zeta_n(x_{i}))-\sum_{i=1}^{\mathcal{N}}y_{i}\zeta_n(x_{i})\right)\right)\Bigg|.
	\end{align*}
	Moreover, by H\"older's inequality ($p=1,\;q=\infty$) and \eqref{2.4.1} we get
	\begin{align*}
	& 	\Bigg|\sum_{m=1}^{N^{\text{train}}}\alpha_m\left(\sum_{n=1}^{\mathcal{R}}\pi_{\zeta_\mathcal{R}}^{n}\left[\phi_{k,m}\right]\left(\sum_{i=1}^{\mathcal{N}}w_{i}\zeta_n(x_{i}))-\sum_{i=1}^{\mathcal{N}}y_{i}\zeta_n(x_{i})\right)\right)\Bigg|
	\\
	\le&\left(\sum_{m=1}^{N^{\text{train}}}|\alpha_m|\right)\max_{m\in\mathbb{M}}\left(\sum_{n=1}^{\mathcal{R}}\Bigg|\pi_{\zeta_\mathcal{R}}^{n}\left[\phi_{k,m}\right]\left(\sum_{i=1}^{\mathcal{N}}w_{i}\zeta_n(x_{i}))-\sum_{i=1}^{\mathcal{N}}y_{i}\zeta_n(x_{i})\right)\Bigg|\right)
	\\
	\le&\epsilon_1\max_{m\in\mathbb{M}}\left(\sum_{n=1}^{\mathcal{R}}|\pi_{\zeta_\mathcal{R}}^{n}\left[\phi_{k,m}\right]|\right)\sum_{m=1}^{N^{\text{train}}}|\alpha_m|,
	\end{align*}
	which implies that 
	\begin{equation}\label{eq8}
	{\bf B} \leq \epsilon_1\max_{m\in\mathbb{M}}\left(\sum_{n=1}^{\mathcal{R}}|\pi_{\zeta_\mathcal{R}}^{n}\left[\phi_{k,m}\right]|\right)\sum_{m=1}^{N^{\text{train}}}|\alpha_m|.	
	\end{equation}
    In conclusion, by using H\"older's inequality ($p,q=2$) for vectors in $\mathbb{R}^{\mathcal{K}}$ and in $\mathbb{R}^{N^{\text{train}}}$, respectively, we obtain 
    \begin{align*}
    {\bf C} \le&\sum_{m=1}^{N^{\text{train}}}|\alpha_m|\Bigg|\sum_{i=1}^{\mathcal{K}}\ws_i\{\phi_{k,m}(\xs_{i})-\Pi_{\zeta_\mathcal{R}}\left[\phi_{k,m}\right](\xs_{i})\}\Bigg|
    \\
    \le&\sum_{m=1}^{N^{\text{train}}}|\alpha_m|\|\ws\|_{2,\mathbb{R}^{\mathcal{K}}}\|\phi_{k,m}-\Pi_{\zeta_\mathcal{R}}\left[\phi_{k,m}\right]\|_{2,\mathbb{R}^{\mathcal{K}}}
    \\
    \le&\|\ws\|_{2,\mathbb{R}^{\mathcal{K}}}\sqrt{\sum_{m=1}^{N^{\text{train}}}\alpha_m^2}\sqrt{\sum_{m=1}^{N^{\text{train}}}\|\phi_{k,m}-\Pi_{\zeta_\mathcal{R}}\left[\phi_{k,m}\right]\|_{2,\mathbb{R}^{\mathcal{K}}}^2}
    \end{align*}
    and, in view of \eqref{eq9} we get
    \begin{align}\label{eq3.5}
    \begin{split}
    {\bf C} \le&\|\ws\|_{2,\mathbb{R}^{\mathcal{K}}}\sqrt{\sum_{m=1}^{N^{\text{train}}}\alpha_m^2}\sqrt{\sum_{m=1}^{N^{\text{train}}}\|\phi_{k,m}-\Pi_{\zeta_\mathcal{R}}\left[\phi_{k,m}\right]\|_{2,\mathbb{R}^{\mathcal{K}}}^2}
    \\ 	
    &\le\|\ws\|_{2,\mathbb{R}^{\mathcal{K}}}\sqrt{\sum_{m=1}^{N^{\text{train}}}\alpha_m^2}\sqrt{\sum_{i=\mathcal{R}+1}^{K\cdot N^{\text{train}}+1}\sigma_i^2}.
    \end{split}
    \end{align}
    Therefore, combining together (\ref{eq12bis}), (\ref{eq8}) and (\ref{eq3.5}) we get (\ref{eqLemma}).
\end{proof}
Some comments on the terms of (\ref{eqLemma}) are required. We have two contributions of different nature: the first one is due to the use of truncated SVD to compress the data while the second one depends on how well the training dataset represents the structure of the parametric manifold. Regarding the first contribution, we also observe that it consists of two terms (\ref{eq12}) and (\ref{eq3.5}), even if we expect one to be negligible with respect to the other. Indeed, we employed (\ref{eq9}) for both the terms, but $E_{SVD}^{\mathcal{R}}$ is the sum of the square $\ell^2$-norm in $\mathbb{R^{\mathcal{N}}}$ while in (\ref{eq3.5}) we actually deal with $\ell^2$-norm in $\mathbb{R^{\mathcal{K}}}$. By definition $\xs_k=x_{i_k}$ with $1\le k\le \mathcal{K}$, therefore we can verify that
\begin{equation}
\|\phi_{k,m}-\Pi_{\zeta_\mathcal{R}}\left[\phi_{k,m}\right]\|_{2,\mathbb{R}^{\mathcal{K}}}\le\|\phi_{k,m}-\Pi_{\zeta_\mathcal{R}}\left[\phi_{k,m}\right]\|_{2,\mathbb{R}^{\mathcal{N}}},\;\;\;\forall k\in\mathbb{K},\;\forall m\in\mathbb{M},
\end{equation}
and moreover, since $K\ll\mathcal{N}$, remembering that the definition of $\ell^2$-norm requires the sum of the squared values in each entry, we expect
\begin{equation}
\|\phi_{k,m}-\Pi_{\zeta_\mathcal{R}}\left[\phi_{k,m}\right]\|_{2,\mathbb{R}^{\mathcal{K}}}\ll\|\phi_{k,m}-\Pi_{\zeta_\mathcal{R}}\left[\phi_{k,m}\right]\|_{2,\mathbb{R}^{\mathcal{N}}},\;\;\;\forall k\in\mathbb{K},\;\forall m\in\mathbb{M};
\end{equation}
which makes contribution of (\ref{eq3.5}) to be expected negligible with respect to (\ref{eq12bis}).\\
Clearly Lemma \ref{lemma1} is not applicable as a priori error estimate, it is necessary to choose an interpolation system. In doing so we demonstrate the following theorem.
\begin{theorem}
	Let
	\begin{align}
	&\Delta\equiv\max_{\mu\in\mathcal{D}}\left(\min_{m\in\mathbb{M}}\|\mu-\mu_{m}^{\text{train}}\|_{2}\right),\label{teo1}\\
	&S_f=\max_{m\in\mathbb{M}}\left(\max_{k\in\mathbb{K}}\left(\sum_{n=1}^{R}\Bigg|<\phi_{k,m},\zeta_{n}>\Bigg|\right)\right).\label{teo3}
	\end{align}
	Suppose the set of functions (\ref{I.1}) satisfies a global Lipschitz condition with respect to the parameters, that is
	\begin{equation}
	\sup_{k\in\mathbb{K}}\sup_{\mu',\mu''\in\mathcal{D}}||f_k(\;\cdot\;,\mu')-f_k(\;\cdot\;,\mu'')||_{L^{\infty}(\Omega)}\le L_f||\mu'-\mu''||_2,
	\label{teo2}
	\end{equation}
	with $L_f$ a positive constant. Then for any $\mu\in\mathcal{D}$, we have that
	\begin{equation}\label{eq:fest}
	\max_{k\in\mathbb{K}}|I_k^{\text{full}}(\mu)-I_k^{sparse}(\mu)|\le\left(\|w\|_{2}+\|\ws\|_{2}\right)\sqrt{\sum_{i=\mathcal{R}+1}^{K\cdot N^{\text{train}}+1}\sigma_i^2}+\epsilon_1 S_f+2|\Omega|L_f\Delta.
	\end{equation}
\end{theorem}
\begin{proof}
	In (\ref{eqLemma}), assume an $\alpha\in\mathbb{R}^{N^{\text{train}}}$ which is sub-optimal
	\begin{equation}
	\alpha_m=\begin{cases}
	\alpha_m=1&\mbox{if }m=\tilde{m}\\
	\alpha_m=0&\mbox{if }m\in\mathbb{M}\setminus \tilde{m}
	\end{cases}
	\label{eqp1}
	\end{equation}
	for $\tilde{m}=\argmin_{m\in\mathbb{M}}\|\mu-\mu_{m}^{\text{train}}\|_2$. From this choice, Lemma \ref{lemma1}, (\ref{teo1}) and (\ref{teo2}) we get the first and third term on the right-hand side of (\ref{eq:fest}). For the second term, starting from (\ref{eqLemma}) by using (\ref{eqp1}), (\ref{2.1.1}), (\ref{teo3}) we obtain
	\begin{align}
	\max_{k\in\mathbb{K}}&\left(\epsilon_1\max_{m\in\mathbb{M}}\left(\sum_{n=1}^{\mathcal{R}}|\pi_{\zeta_\mathcal{R}}^{n}\left[\phi_{k,m}\right]|\right)\sum_{m=1}^{N^{\text{train}}}|\alpha_m|\right)=\epsilon_1\max_{k\in\mathbb{K}}\left(\left(\sum_{n=1}^{\mathcal{R}}|\pi_{\zeta_\mathcal{R}}^{n}\left[\phi_{k,\tilde{m}}\right]|\right)\right)\nonumber\\
	&\le\epsilon_1\max_{k\in\mathbb{K}}\left(\max_{m\in\mathbb{M}}\left(\sum_{n=1}^{\mathcal{R}}\Bigg|<\phi_{k,m},\zeta_{n}>\Bigg|\right)\right)=\epsilon_1 S_f.
	\end{align}
\end{proof}
As observed in \cite{p1}, it is natural to add the hypothesis $\Delta\rightarrow0$ as $N^{\text{train}}\rightarrow\infty$. This implies that for a large training dataset we expect the contribution of third term in (\ref{eq:fest}) to be negligible. The $\ell^2$-norm of the sparse rule, that is only available a posteriori, can be bounded by its $\ell^1$-norm that is equal to $|\Omega|$. The constant $S_f$ can be numerically computed once $\mathcal{R}$ is fixed, indeed
\begin{equation}\label{eq:appSf}
	S_f\approx \max_{j}\left\{\sum_{i=1}^{\mathcal{R}}\Big|\left(U(:, 1:\mathcal{R})^TA^T\right)_{i,j}\Big|\right\}.
\end{equation}
Therefore, if one wants to satisfy (\ref{eq1.8}) it is enough to choose, for instance, $\mathcal{R}$ and $\epsilon_1$ such that the sum of the first and second term of (\ref{eq:fest}) is equal or less than $\varepsilon/2$. Since we have the sum of two terms one can play in making one negligible with respect to the other. This could promote a smaller value of $\mathcal{R}$ or a larger norm of the residual that enforces sparsity. We discuss this argument with the help of numerical tests in the next section.
\section{Numerical Examples}
We consider two examples: the first is related to the evaluation in space and time of the fundamental solution of the one dimensional linear Schrödinger equation with Gaussian initial data; the second concerns the reduced-basis method. All the computations are performed on a laptop with $2.60$ GHz Intel Core i7 processor, using \verb|Matlab 2020a|, also the associated codes are available in \cite{p19}.
\subsection{Numerical approximation of the 1D Schrödinger fundamental solution}
\begin{figure}[t]
		
		\centering{
			\subfigure{
				\includegraphics[width=0.48\textwidth]{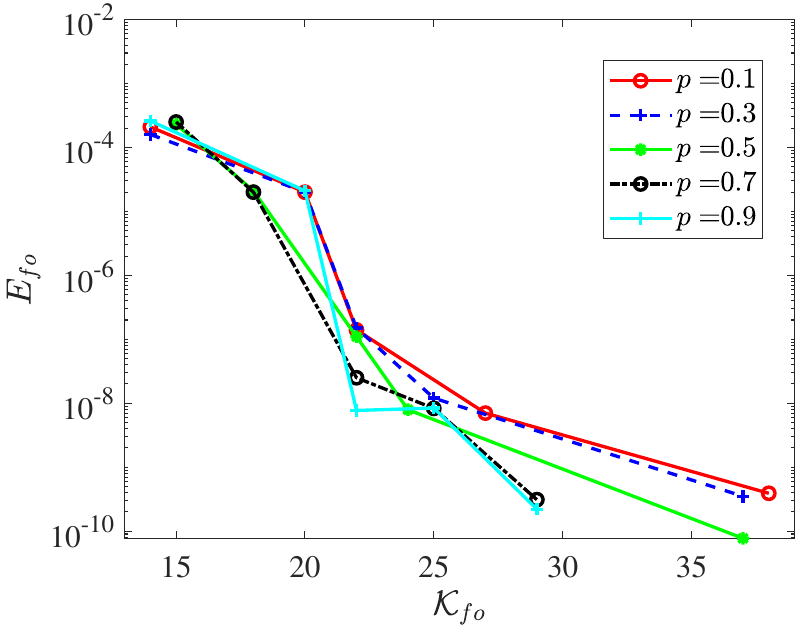}}
			%\hspace{0.05cm}
			\subfigure{
				\includegraphics[width=0.48\textwidth]{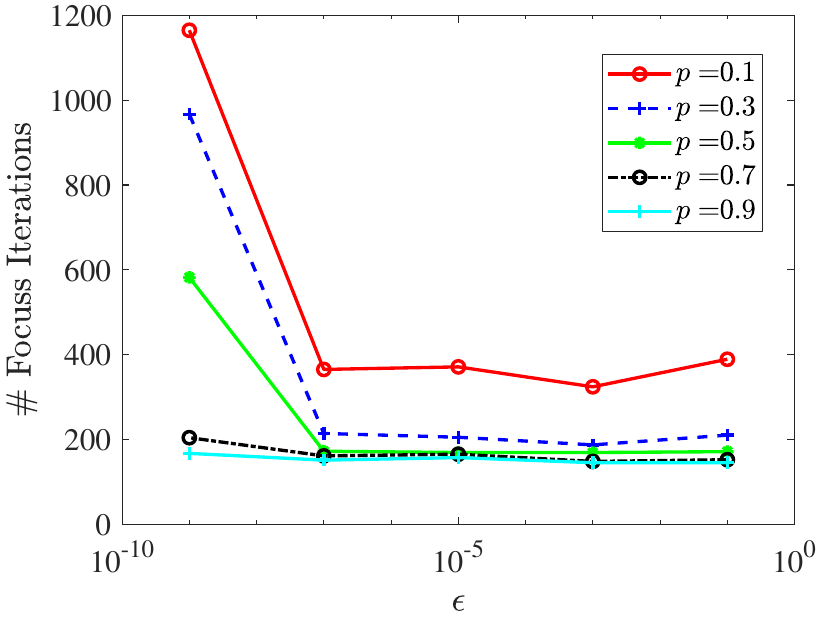}}
		}
		
		\caption{Schrödinger fundamental solution test problem: error (\ref{eq:err}) with respect to the number of quadrature nodes (left) and total \texttt{FOCUSS} iterations with respect to the target accuracy $\epsilon$ (right) for different values of $p$.}
		\label{fig2}
\end{figure}

\begin{figure}[]

	\centering{
		\subfigure{
			\includegraphics[width=0.48\textwidth]{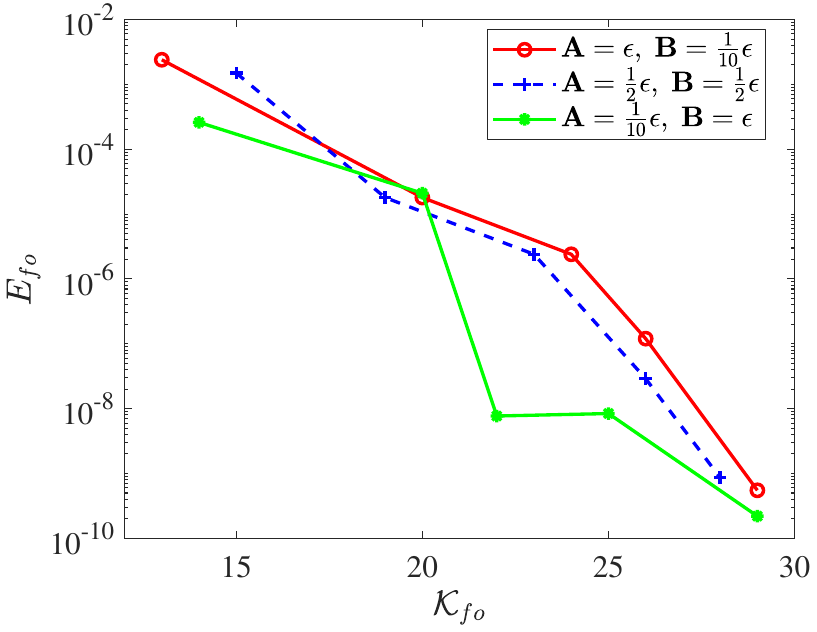}}
		%\hspace{0.05cm}
		\subfigure{
			\includegraphics[width=0.48\textwidth]{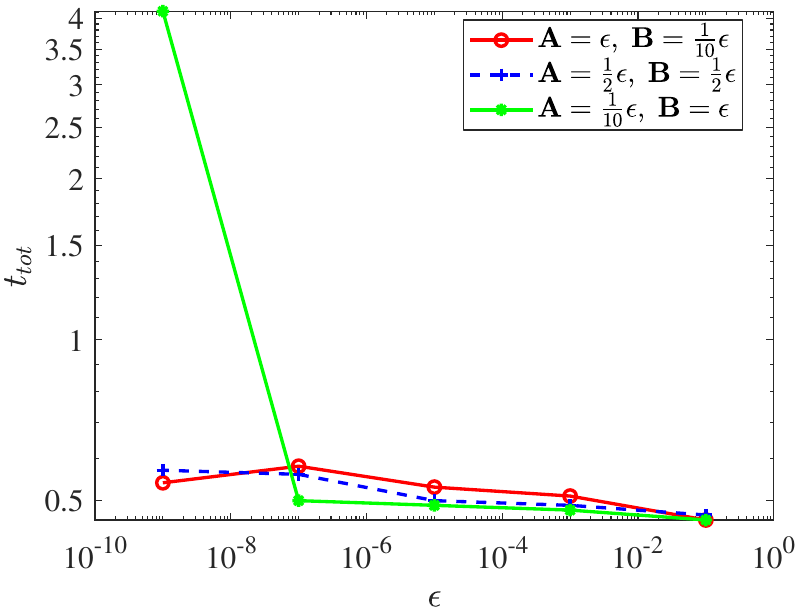}}
	}
	
	\caption{Schrödinger fundamental solution test problem: error (\ref{eq:err}) with respect to the number of quadrature nodes (left) and offline computational time with respect to target accuracy (right) for the three strategies to recover quadrature rule of order of accuracy $\epsilon$.}
	\label{fig3}
\end{figure}

We consider the Cauchy problem associated with the $1$-dimensional linear Schrödinger equation
\begin{equation} \label{eqSc}
\begin{cases}
     i\partial_t \psi(x,t)=-\partial_{xx} \psi(x,t),       \\
	\psi(x,0)=\psi_0(x).\\
\end{cases}
\end{equation}
The so-called fundamental solution of (\ref{eqSc}) is given by
\begin{equation}\label{eqFu}
 \psi(x,t)=\left(\frac{1}{\sqrt{4\pi i t}}\right) e^{\frac{i|x|^2}{4t}}\int_{\mathbb{R}}e^{-\frac{ixy}{2t}}e^{\frac{i|y|^2}{4t}}\psi_0(y)\;\text{d}y,
\end{equation}
that, up to rescaling and to multiplication by a function of modulus $1$, is the Fourier transform of the initial condition.\\
Setting $K=1$ and $\mu\equiv(\mu_{1},\mu_{2})\equiv(x,t)\in\mathcal{D}\subset\mathbb{R}^2$, we have that (\ref{eqFu}) can be view as a parametrized integral, i.e. $\psi(x,t)=I(\mu)$ of (\ref{eq1.2}) with integrand function
\begin{equation}\label{eqIn}
f(y,\mu)\equiv e^{-\frac{ixy}{2t}}e^{\frac{i|y|^2}{4t}}\psi_0(y).
\end{equation} 
We take the initial data $\psi_0(y)=e^{-\frac{y^2}{2}}$, so that the real part $\Re({f(y,\mu)})$ is symmetric while the imaginary part $\Im({f(y,\mu)})$ is antisymmetric. Therefore the integrand function (\ref{eqIn}) becomes
\begin{equation}\label{eqIn2}
	f(y,\mu)\equiv \left(\cos\left({-\frac{xy}{2t}}\right)\cos\left({\frac{|y|^2}{4t}}\right)-\sin\left({-\frac{xy}{2t}}\right)\sin\left({\frac{|y|^2}{4t}}\right)\right)e^{-\frac{y^2}{2}},
\end{equation}
the domain of integration is reduced to $\Omega=[0, \infty)$, and we finally deal with
\begin{equation}\label{eqFu2}
	I(\mu)=\int_{0}^{\infty}f(y,\mu)\;\text{d}y.
\end{equation}
A quadrature formula is required to numerically treat (\ref{eqFu2}); here we make use of a trapezoidal rule over $\mathcal{N}$ equally spaced points on the interval $[0, y_{max}]$. Note that the error between (\ref{eqFu2}) and the approximated integral $I^{\text{full}}(\mu)$ has also a component due to the truncation at $y_{max}$ of the integration domain. However this component can be arbitrary small since we have an exponential decay to zero of  $f(x,\mu)$ as $x\longrightarrow\infty$.\\
We set $t\in[0,4]$ and $x\in[0.2,2]$ such that $\mathcal{P}\equiv[0,4]\times[0.2,2]$. The full quadrature is applied with $y_{max}=4$ and $\mathcal{N}=1200$. We consider the training dataset $\Xi^{\text{train}}$ from uniform sampling over $\mathcal{D}$ of $N^{\text{train}}=J\times J$ elements. The sparse quadrature rule $\{\hat{y},\ws\}$ is obtained as post processing of our routine based on \verb|FOCUSS|, the error it introduces is measured as
\begin{equation}\label{eq:err}
E(\Xi^{\text{test}})=\max_{\mu\in\Xi^{\text{test}}}|I^{\text{full}}(\mu)-I^{\text{sparse}}(\mu)|.
\end{equation}
where $\Xi^{\text{test}}$ is the parameter test sample of size $200^2$ constructed as the tensorization of a uniformly random distributed grid of size $200$ in each of the two parameter directions. We consider $J=40$ for the construction of the train sample.\\ 
First of all we analyse the influence of $p$ in the recover of the sparse quadrature rule and in the convergence of the \texttt{FOCUSS} algorithm. Figure \ref{fig2} (left) shows the error (\ref{eq:err}) over the number of quadrature points for values of $p$ between $0$ and $1$. We do not see remarkable difference except for the smallest error where the highest values of $p$ require less nodes to reach the target accuracy. Figure \ref{fig2} (right) confirms that the choice of highest values of $p$ is ``optimal'' since it also minimizes the number of iterations for the algorithm convergence.
\\
Next we discuss about the choice of $\mathcal{R}$ and $\epsilon_1$ to get a sparse quadrature rule with integration error (\ref{eq:err}) smaller then $\epsilon$. Let us assume to have an enough large train sample, so that the third component of (\ref{eq:fest}) is much smaller than $\epsilon$. We would like to have
\begin{equation}
	E(\Xi^{\text{test}})\le \epsilon=\underbrace{\Big(\|w\|_{2}+|\Omega|\Big)\sqrt{\sum_{i=\mathcal{R}+1}^{K\cdot N^{\text{train}}+1}\sigma_i^2}}_{\mathbf{A}}+\underbrace{\epsilon_1 S_f}_{\mathbf{B}},
\end{equation}
where $S_f$ is approximated according (\ref{eq:appSf}). Once the singular value decomposition of the train matrix is computed, there are three possible strategies that can be exploited, concerning the choice of $\mathcal{R}$ and $\epsilon_1$:
\begin{enumerate}
\item $\mathbf{A}=\rm{O}(\epsilon)$ and $\mathbf{B}\ll\epsilon$;
\item $\mathbf{A}=\frac{1}{2}\epsilon$ and $\mathbf{B}=\frac{1}{2}\epsilon$;
\item$\mathbf{A}\ll\epsilon$ and $\mathbf{B}=\text{O}(\epsilon)$.
\end{enumerate}
To implement first and third option we chose to make $\mathbf{B}$ and $\mathbf{A}$ ten times smaller than the target accuracy. Figure \ref{fig3} (left) shows the results in terms of error over sparsity for the three strategies listed. We report that, while all the three options are effective in recovering sparse rules of the required tolerance, the third one recovers sparse rules which provide smaller error with a similar number of quadrature points. The rightmost panel in Figure \ref{fig3} shows that the computational time to run the three strategies is similar, expect that for the third one for which it is considerably higher for $\epsilon=10^{-9}$. Based on these results, we use the second option when $\epsilon\le10^{-7}$ and the third one for all the other values of $\epsilon$. 
\begin{figure}[]

	\centering{
		\subfigure{
			\includegraphics[width=0.48\textwidth]{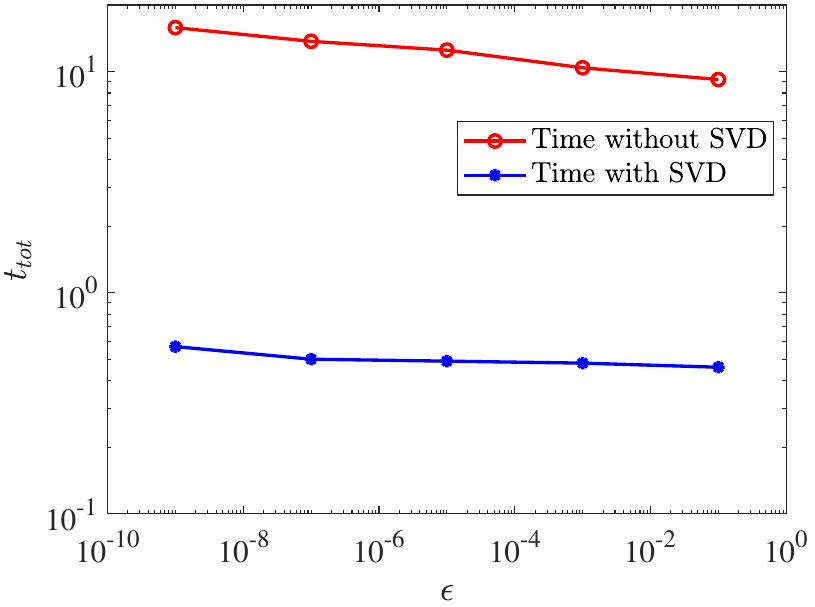}}
		%\hspace{0.05cm}
		\subfigure{
			\includegraphics[width=0.48\textwidth]{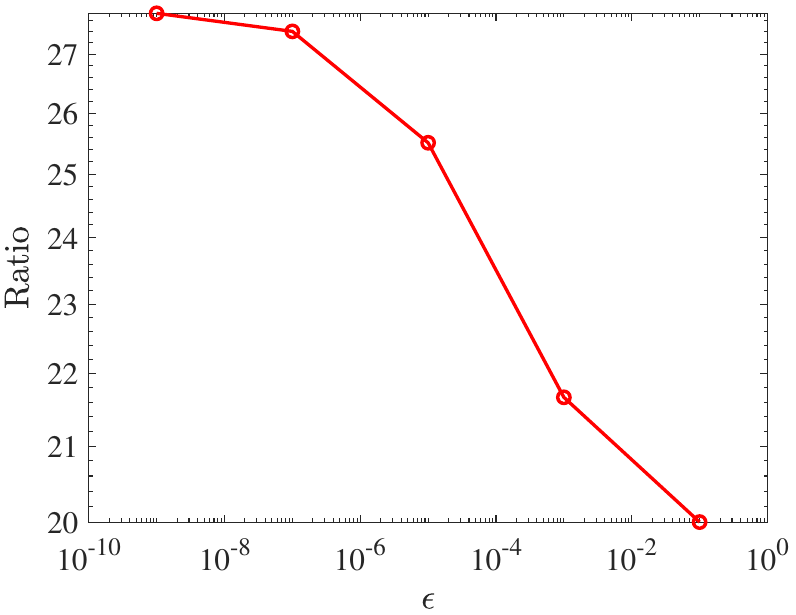}}
	}
	
	\caption{Schrödinger fundamental solution test problem: computational time to recover the sparse rule with respect to the target accuracy. Direct comparison between the strategies (left) and ratio between the recovering time (right) without SVD and with SVD.}
	\label{fig4}
\end{figure}\\
Figure \ref{fig4} shows the advantage of performing a truncated SVD on the constrain matrix, indeed our strategy results to be between the 20 and 29 times faster than a direct application of the \texttt{FOCUSS} algorithm without truncated SVD.
We set different values of the integration tolerance $\epsilon$ and we report in Table \ref{T1} ($J=40$) the results collected. In the same table we also compare our strategy with those methods which are able to solve problem (\ref{eq1.5}) under constrains (\ref{eq1.6}-\ref{eq1.6bis}), these are: the dual-simplex algorithm \cite{p10} of the linear programming routine of \texttt{Matlab} and the empirical cubature method \cite{p18} that was also implemented in \texttt{Matlab}.\\
The comparison among the different sparse quadrature recovery strategies is carried out in terms of the following criteria, which are listed in order of importance:
\begin{enumerate}
	\item The ability of the algorithm to strictly enforce the target accuracy tolerance imposed to the sparse rule. Indeed, the a priori estimation of the error is a key feature in the applications mentioned in the scope of this work;
	\item The ability to produce the sparsest quadrature rule for a fixed accuracy. Given two methods which provide sparse rules at the same target accuracy, we prioritize the one that gives rules with fewer nodes since this will speed-up the online phase in the applications;
	\item The overall execution time for the sparse quadrature recovery. Methods that are faster in providing the sparse rule of the required accuracy (in the offline phase) are to be preferred;
	\item Finally, given two sparse rules with same number of quadrature nodes originated by two different methods and both respecting the first criterion, we promote the one which results in a smaller integration error (\ref{eq:err}).
\end{enumerate}
The results of our numerical experiments are presented in Table \ref{T1}. They show that our strategy is always able to recover sparse rules of the required accuracy. The linear programming strategy succeeds only for $\epsilon=10^{-1}$, while the empirical cubature is effective until $\epsilon=10^{-5}$ and it is not convergent when $\epsilon=10^{-9}$. Figure \ref{fig8} compares the algorithms in terms of error (\ref{eq:err}) with respect to sparsity. Our strategy based on \texttt{FOCUSS} algorithm always recovers a more accurate quadrature rule for a fixed number of quadrature points. We also observed that the empirical cubature was not able to return a sparse rule when high accuracies were demanded, we suspect that this is due to the limited accuracy reachable by the \texttt{lsqnonneg} function implemented by \texttt{Matlab}.  
\begin{table}[t]
	\centering
	\begin{tabular}{c|ccccc}
		
		$\epsilon$&$1\cdot10^{-1}$& $1\cdot10^{-3}$ & $1\cdot10^{-5}$ & $1\cdot10^{-7}$& $1\cdot10^{-9}$ \\
		\hline
		$\mathcal{K}_{lp}$ & $9$    & $13$ & $18$ & $22$& $27$\\
		
		$\mathcal{K}_{ec}$ & $9$    & $16$ & $20$ & $21$& NA\\
		
		$\mathcal{K}_{fo}$    &$15$ & $18$& $22$ & $24$&$28$\\
		
		$E_{lp}$ & $1.0\cdot10^{-1}$& $2.0\cdot10^{-3}$& $1.5\cdot10^{-4}$ & $7.3\cdot10^{-6}$&$3.8\cdot10^{-9}$\\
		
		$E_{ec}$ & $2.0\cdot10^{-2}$& $5.3\cdot10^{-4}$& $1.3\cdot10^{-5}$ & $1.2\cdot10^{-3}$&NA\\
		
		$E_{fo}$  &$2.6\cdot10^{-4}$& $2.1\cdot10^{-5}$ & $7.8\cdot10^{-9}$ &$8.4\cdot10^{-9}$&$8.6\cdot10^{-10}$\\
		
	\end{tabular}
	\caption{Schrödinger fundamental solution test problem: results in terms of number of quadrature nodes and errors (\ref{eq:err}) for linear programming ($lp$), empirical cubature ($ec$) and \texttt{FOCUSS} based strategy ($fo$).}
	\label{T1}
\end{table}
\begin{figure}[]

	\centering{
		
		\includegraphics[width=0.48\textwidth]{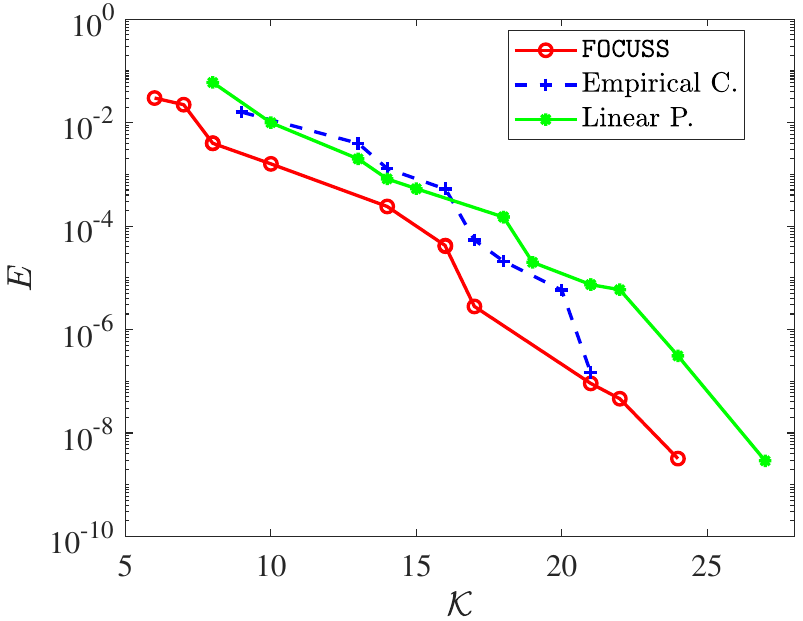}
	}
	
	\caption{Schrödinger fundamental solution test problem: error (\ref{eq:err}) with respect to number of quadrature points for the three sparse recovery algorithms.}
	\label{fig8}
\end{figure}\\
Finally, we report the average computational time required by each method to find the sparse rule for different sizes of the training dataset $|\Xi^{\text{train}}|=J\times J$. Results are displayed in Table \ref{T4}. The fastest method is the empirical cubature except for $\epsilon=10^{-9}$, while the linear programming is the lowest in any case tested. We observe that the algorithm proposed and the empirical cubature require more time as $\epsilon$ is lowered while the one implementing the $\ell^1$-norm minimization has an almost constant execution time with respect to $\epsilon$. Moreover, the empirical cubature fails to converge for $\epsilon= 10^{-9}$ and $J\ge40$. Another feature displayed concerns the comparison of the methods in terms of increasing computational time as the size of the problem becomes bigger. Indeed, the results suggest that the algorithms based on \verb|FOCUSS| and the empirical cubature have a better scaling in this sense, compared to the one implementing the $\ell^1$-norm minimization. We also note that most of the computational time demanded by our strategy is due to the initial truncated SVD which cost does not scale linearly with respect to $J$, anyway for very large datasets we can take advantage of randomized algorithms for matrix decomposition, see \cite{p20}, to strongly speed up the computation.
\begin{table}[t]
	\centering
	\begin{tabular}{|c|ccc|ccc|ccc|}
		\hline
		\multirow{2}{*}{$\epsilon$}&\multicolumn{3}{c|}{$J=20$}&\multicolumn{3}{c|}{$J=40$}&\multicolumn{3}{c|}{$J=80$}\\
		\cline{2-10}
		& $t_{lp}$ &$t_{ec}$& $t_{fo}$ & $t_{lp}$ &$t_{ec}$& $t_{fo}$& $t_{lp}$ &$t_{ec}$& $t_{fo}$\\
		\hline
		\hline
		$10^{-1}$   & 0.49    &0.02&0.25&
		3.07&0.04&0.67&
		38.7&0.23&2.8 \\
		
		$10^{-3}$   & 0.47    &0.04&0.34&
		3.23&0.09&0.78&
		45.6&0.49&2.9 \\
		
		$10^{-5}$    & 0.50 &0.03&0.34&
		3.01&0.09&0.69&
		37.9&0.51&3.03\\
		
		$10^{-7}$    & 0.77 &0.09&0.41&
		3.14&0.23&0.68&
		39.6&0.64&2.90\\
		
		$10^{-9}$  & 0.77 &5.65&0.46&
		4.89&$\infty$&1.05&
		61.2&$\infty$&3.84\\
		
		\hline
	\end{tabular}
	\cprotect\caption{Schrödinger fundamental solution test problem: average time(s) to compute the sparse rule with linear programming ($lp$), empirical cubature ($ec$) and \verb|FOCUSS| based method ($fo$) for different sizes of the training dataset $|\Xi^{\text{train}}|=J\times J$ and different integration tolerance $\epsilon$.}
	\label{T4}
\end{table}
\subsection{Nonlinear reduced-basis diffusion example}
The second example is related to the framework of reduced-basis method \cite{p11}. We consider the following parametrized nonlinear diffusion problem: for a given $\mu=(\mu_{1},\mu_{2},\mu_{3})\in\mathcal{D}\subset\mathbb{R}^3$, find $u(x,\mu)\in\mathcal{V}\equiv H^1_0(\Omega)\equiv\{v\in H^1(\Omega)|v|_{\partial \Omega}=0\}$ such that 
\begin{equation}
	r(u(\mu),v;\mu)=\int_{\Omega}\kappa(u,\mu)\nabla u\cdot\nabla v \;\text{d}x-\int_{\Omega}10v\;\text{d}x=0\;\;\forall v\in\mathcal{V};
\end{equation}
where $\Omega=[0,1]^2$ is the domain and $\kappa$ is a nonlinear diffusion coefficient defined as: 
\begin{equation}
\kappa(u,\mu)=\begin{cases} 1+\mu_1u, & \mbox{if }x\in\Omega_1 \\ \mu_2+\mu_3u^2, & \mbox{if }x\in \Omega\setminus\Omega_1
\end{cases}
\end{equation}
with $\Omega_1$ is a circle of radius $0.25$ centred in $\Omega$ and $\mu_1\in[0,10]$, $\mu_2\in[1,10]$, $\mu_{3}\in[0,10]$. The space $\mathcal{V}$ is equipped with the standard $H^1_0$ inner product and norm.\\
To solve this problem we use a standard linear finite-element discretization, i.e we introduce the space $\mathcal{V}_h\equiv\{v\in\mathcal{V}|v|_k\in\mathbb{P}^1(k),\forall k\in\mathcal{T}_h\}\subset\mathcal{V}$ where $\mathcal{T}_h$ is a triangulation over $\Omega$ of $N_{ele}=8288$ elements, see Figure \ref{fig1} (left). In this way we can solve the finite dimensional problems: given $\mu\in\mathcal{D}$, find $u_h(\mu)\in\mathcal{V}_h$ such that
\begin{equation}
	r(u_h(\mu),v;\mu)=0\;\;\forall v\in\mathcal{V}_h;
	\label{eq5.1}
\end{equation}
which integrals are evaluated by a full quadrature rule that consists of $\mathcal{N}=33152$ quadrature points. For fixed values of $\mu_1$, $\mu_2$ and $\mu_3$ we solve the nonlinearity by means of Picard iterations. \\
As next step we introduce a reduced-basis approximation of problem (\ref{eq5.1}) by defining the $\mathcal{V}_{N}\equiv\text{span}\{\zeta_i\}_{i=1}^{N}\subset\mathcal{V}_h$ and then we can state the reduced-basis (full quadrature) version of our problem: given $\mu\in\mathcal{D}$, find $u_N(\mu)\in\mathcal{V}_N$ such that
\begin{equation}
	r(u_N(\mu),v;\mu)=0\;\;\forall v\in\mathcal{V}_N.
	\label{eq5.2}
\end{equation} 
We construct our reduced-basis approximation space by the use of Proper Orthogonal Decomposition (POD). We introduce a training dataset $\Xi^{\text{rb},\text{train}}\subset\mathcal{D}$ of $|\Xi^{\text{rb},\text{train}}|=7^3$ points and we compute and collect the solution of problem (\ref{eq5.1}) $\forall\mu\in\Xi^{\text{rb},\text{train}}$. Then we apply truncated SVD to the data collected in order to extract the most significant modes and define the reduced-basis approximation spaces $\mathcal{V}_{N=1}\subset...\subset\mathcal{V}_{N=N_{max}}$, see for instance \cite{p11}.\\
Finally, we consider a sparse quadrature approximation of (\ref{eq5.2}) and introduce the following residual form:
\begin{equation}
	r^{fo}(u^{fo}_N,v;\mu)=\sum_{i=1}^{\mathcal{K}}\ws_i\kappa(u_N^{{fo}}(\xs_i),\mu)\nabla u^{{fo}}_N(\xs_i)\nabla v(\xs_i)-\int_{\Omega}10v\;\text{d}x=0\;\;\forall v\in\mathcal{V}_N.
\end{equation}
The reduced-basis approximation associated with the reduced quadrature is defined in the following way: given $\mu\in\mathcal{D}$, find $u_N^{{fo}}(\mu)\in\mathcal{V}_N$ such that
\begin{equation}
	r^{hr}(u_N^{{fo}}(\mu),v;\mu)=0\;\;\forall v\in\mathcal{V}_N.
\end{equation} 
We train our sparse quadrature rule requiring that: $\forall\mu\in\Xi^{\text{rb},\text{train}}$ the residual evaluated on the full reduced solution has to be integrated up to an accuracy $\epsilon$, i.e.
\begin{equation}
	|r^{{fo}}(u_N(\mu),v;\mu)|\le\epsilon\;\;\forall v\in \mathcal{V}_{N}.
\end{equation}
This formulation corresponds to $K=N$ and $N^{\text{train}}=|\Xi^{\text{rb},\text{train}}|$, for a total of $N\cdot|\Xi^{\text{rb},\text{train}}|$ constraints to which we must add (\ref{eq1.8bis}). To test the sparse quadrature found with our routine we define a set $\Xi^{\text{test}}$ that consists of $1000$ random uniformly distributed points over $\mathcal{D}$.
\begin{figure}[t]
	
	\centering{
		\subfigure{
			\includegraphics[width=0.48\textwidth]{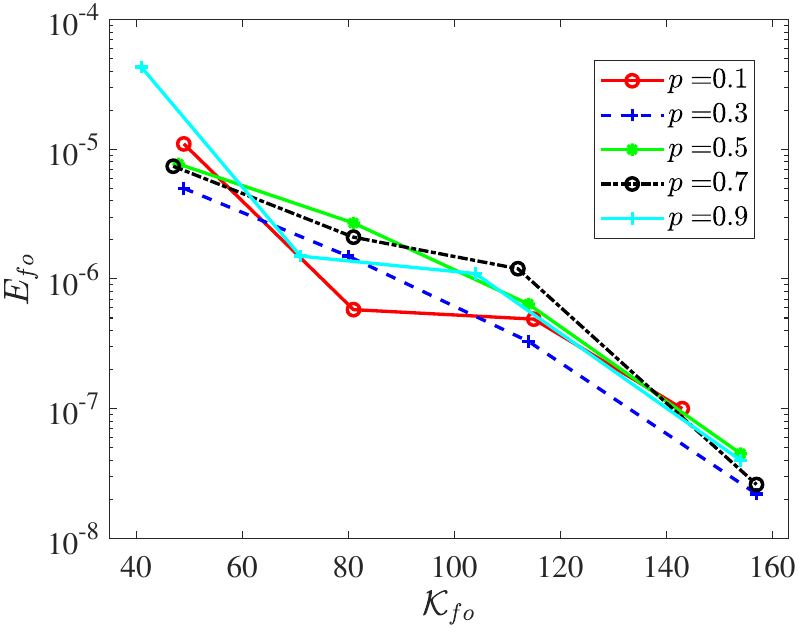}}
		%\hspace{0.05cm}
		\subfigure{
			\includegraphics[width=0.48\textwidth]{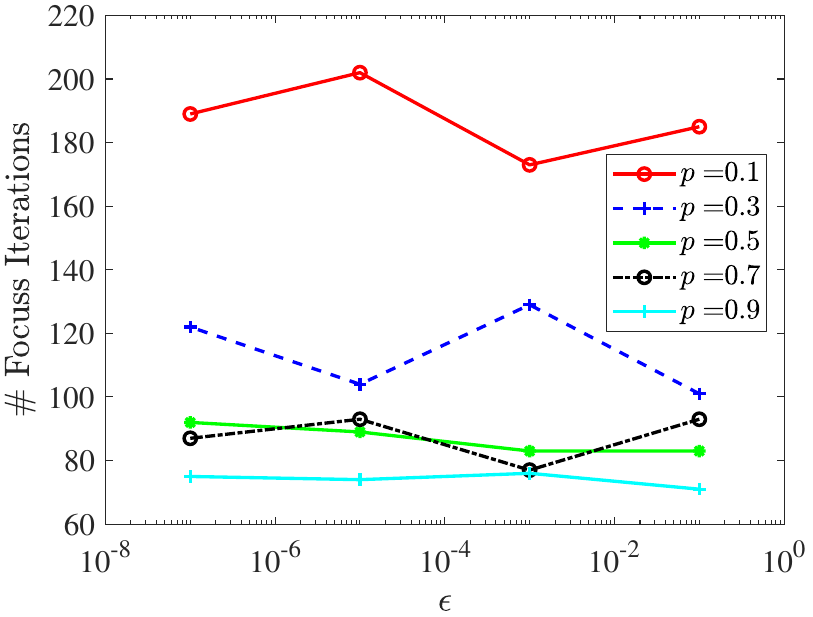}}
	}
	
	\caption{Parametrized nonlinear diffusion test problem: error with respect to the number of quadrature nodes (left) and total \texttt{FOCUSS} iterations with respect to the target accuracy $\epsilon$ (right) for different values of $p$.}
	\label{fig5}
\end{figure}
\begin{figure}[]

	\centering{
		\subfigure{
			\includegraphics[width=0.48\textwidth]{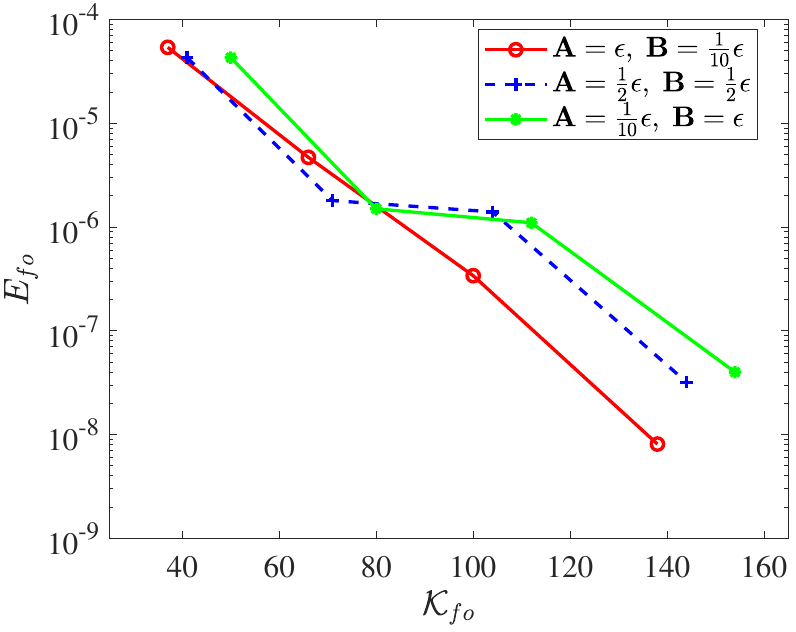}}
		%\hspace{0.05cm}
		\subfigure{
			\includegraphics[width=0.48\textwidth]{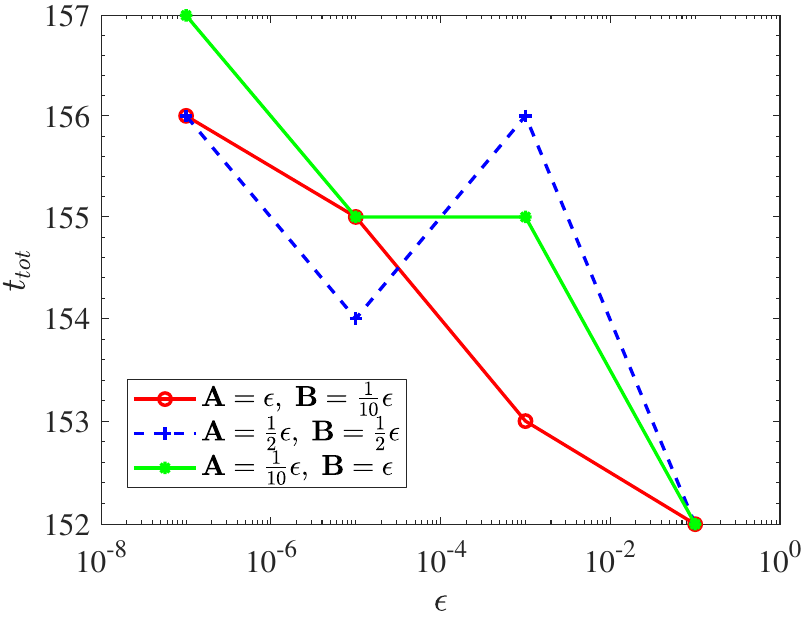}}
	}
	
	\caption{Parametrized nonlinear diffusion test problem: error with respect to the number of quadrature nodes (left) and offline computational time with respect to target accuracy (right) for the three strategies to recover quadrature rule of order of accuracy $\epsilon$.}
	\label{fig6}
\end{figure}\\
We start repeating the numerical experiment done for the Schrödinger fundamental solution example. In Figure \ref{fig6} we analyse the influence of $p$ in recovering sparse quadrature rules and in the convergence of the \texttt{FOCUSS} algorithm. Again the choice of the larger values of $p$ seems to be optimal since the number of iterations for convergence becomes smaller while the behaviour of the error with respect to sparsity does not suffer of remarkable changes. Figure \ref{fig6} (left) shows the results in terms of error over sparsity for the three strategies to integrate with precision $\epsilon$. Again we see that all the three options are effective in recovering sparse rules of the required tolerance. The first option is the one recovering sparse rules which provide smaller errors with also a slightly smaller number of quadrature points. In the right plot of Figure \ref{fig6} we can see that the computational time to run the three strategies is similar. Finally in Figure \ref{fig7} we compare the approach with and without truncated SVD. Again the use of SVD results in a faster algorithm, with a speed between the 6 and 8.5 higher than a direct application of \texttt{FOCUSS} without truncated SVD on the constrain matrix.
\begin{figure}[]

		\centering{
			\subfigure{
				\includegraphics[width=0.48\textwidth]{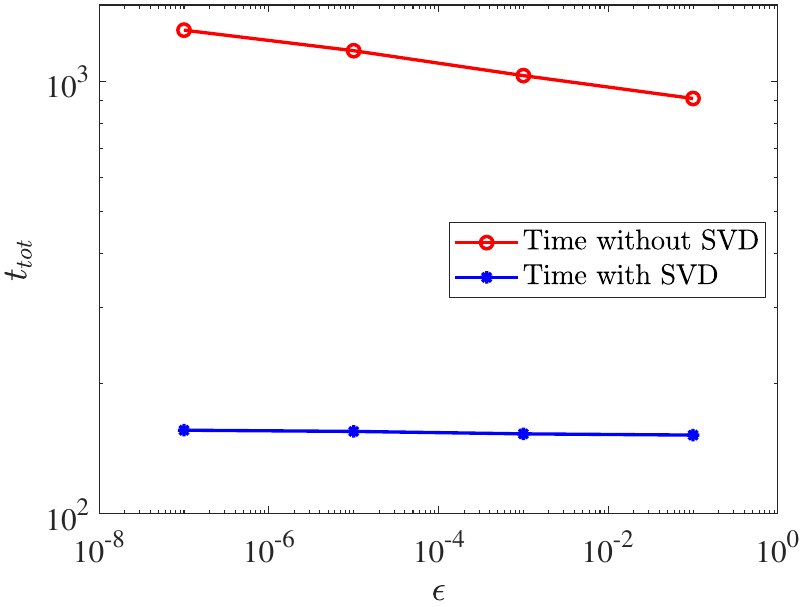}}
			%\hspace{0.05cm}
			\subfigure{
				\includegraphics[width=0.48\textwidth]{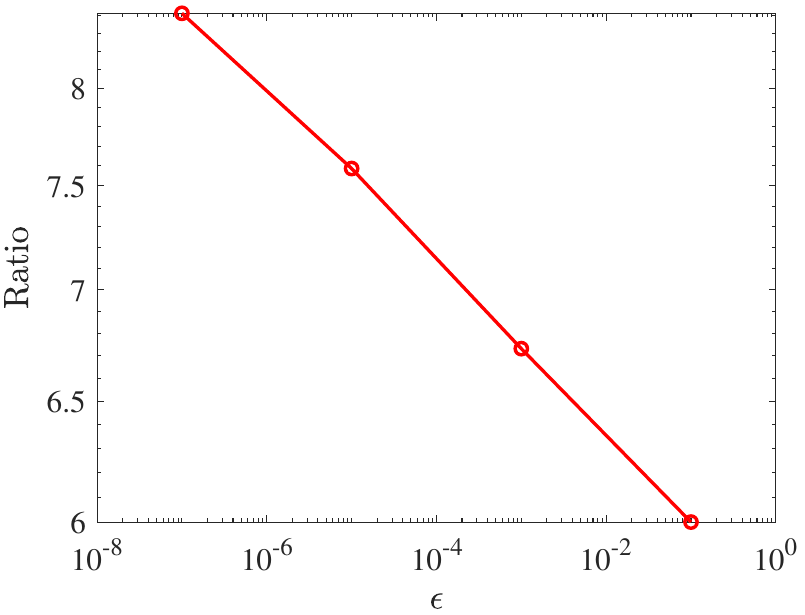}}
		}
		
		\caption{Parametrized nonlinear diffusion test problem: computational time to recover the sparse rule with respect to the target accuracy. Direct comparison between the strategies (left) and ratio between the recovering time (right) without SVD and with SVD.}
		\label{fig7}
\end{figure}\\
In Table \ref{Tab-1} we present the results for the nonlinear diffusion problem using both the reduced-basis approximation and the sparse quadrature rule for $\epsilon=10^{-5}$. The number of reduced quadrature points goes from $5$, when $N=1$, to $136$, when $N=7$, and the quadrature error over the residual is always kept far behind the required precision $\epsilon=10^{-5}$. This reflects in a small error between the reduced-basis approximation using the full quadrature $u_N(\mu)$ and the one using the sparse quadrature $u_N^{fo}(\mu)$. Moreover, this error is always smaller than the error in the reduced-basis approximation (last column of Table \ref{Tab-1}). We observe that $\mathcal{K}$ increases with $N$. This is due to the fact that the number of constraints satisfied by the reduced rule also increases with $N$, implying a rule with an higher number of non zero weights. In Figure \ref{fig1} (right) we can see an example of reduced quadrature nodes distribution along the domain $\Omega$.
\begin{table}[t]
	\begin{center}
		\begin{tabular}{llccc}
			\cline{1-5}
			\multirow{3}{*}{$N$}&\multirow{3}{*}{$\mathcal{K}$}&\multirow{3}{*}{$E(\Xi^{\text{test}})$}&\multirow{3}{*}{$\max\limits_{\substack{\mu\in\Xi^{\text{test}}}}\dfrac{\norm{u_N(\mu)-u_N^{{fo}}(\mu)}_{\mathcal{V}}}{\norm{u_N(\mu)}_{\mathcal{V}}}$}&\multirow{3}{*}{$\max\limits_{\substack{\mu\in\Xi^{\text{test}}}}\dfrac{\norm{u_h(\mu)-u_N^{{fo}}(\mu)}_{\mathcal{V}}}{\norm{u_h(\mu)}_{\mathcal{V}}}$}\\
			&&&&
			\\
			&&&&
			\\
			\cline{1-5}
			&&&&
			\\
			$1$&$5$&$7.1\cdot10^{-7}$&$1.3\cdot10^{-5}$&$2.8\cdot10^{-1}$\\
			
			$2$&$19$&$3.7\cdot10^{-7}$&$2.4\cdot10^{-6}$&$1.6\cdot10^{-1}$\\
			
			$3$&$39$&$2.8\cdot10^{-6}$&$6.3\cdot10^{-6}$&$1.9\cdot10^{-2}$\\
			
			$4$&$64$&$1.2\cdot10^{-6}$&$2.7\cdot10^{-6}$&$7.3\cdot10^{-3}$\\
			
			$5$&$85$&$7.8\cdot10^{-7}$&$1.6\cdot10^{-6}$&$5.0\cdot10^{-3}$\\
			
			$6$&$109$&$6.7\cdot10^{-7}$&$1.2\cdot10^{-6}$&$2.4\cdot10^{-3}$\\
			
			$7$&$136$&$1.2\cdot10^{-6}$&$2.2\cdot10^{-6}$&$4.8\cdot10^{-4}$\\
			
			\cline{1-5}
		\end{tabular}   
	\end{center}
	\caption{Parametrized nonlinear diffusion test problem: the number of reduced-basis functions $N$, the number of reduced quadrature points $\mathcal{K}$, the maximum quadrature error for the nonlinear term over the test set $\Xi^{\text{test}}$, $E(\Xi^{\text{test}})$, the maximum relative difference in the reduced-basis approximation using the full quadrature and sparse quadrature,  $\max_{\mu\in\Xi^{\text{test}}}\dfrac{\norm{u_N(\mu)-u_N^{{fo}}(\mu)}_{\mathcal{V}}}{\norm{u_N(\mu)}_{\mathcal{V}}}$, and the relative error in the reduced-basis approximation using the reduced quadrature, $\max_{\mu\in\Xi^{\text{test}}}\dfrac{\norm{u_h(\mu)-u_N^{{fo}}(\mu)}_{\mathcal{V}}}{\norm{u_h(\mu)}_{\mathcal{V}}}$.}
	\label{Tab-1}
\end{table}
\begin{figure}
	\centering{
		\subfigure{
			\includegraphics[width=0.4\textwidth]{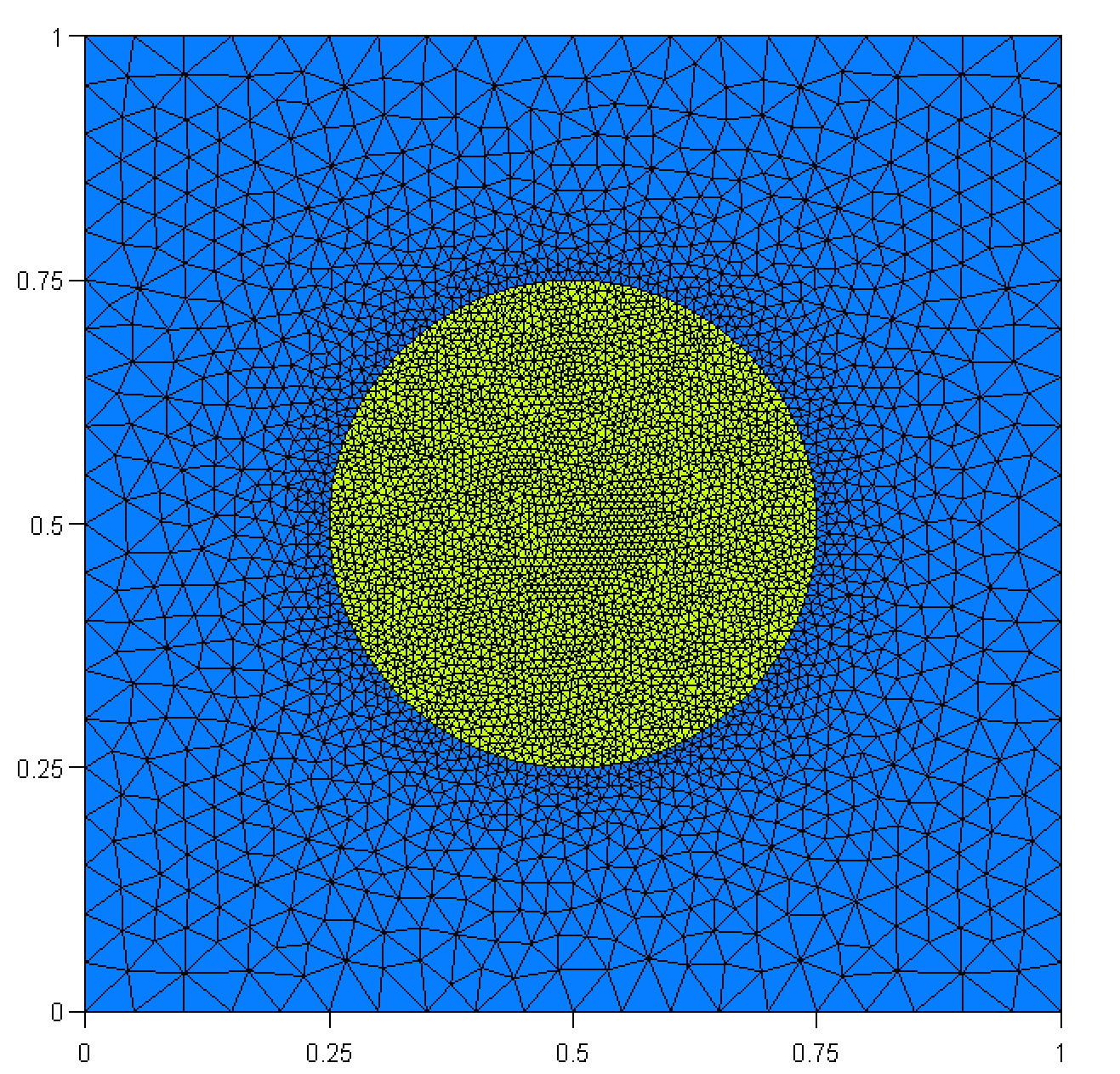}}
		%\hspace{0.05cm}
		\subfigure{
			\includegraphics[width=0.4\textwidth]{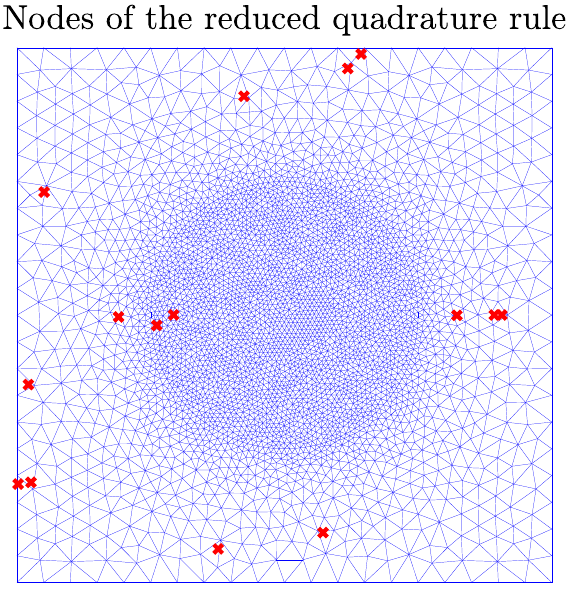}}
	}	
	\caption{Mesh considered for the domain $\Omega$ of the nonlinear diffusion problem (left). Distribution of the reduced quadrature rule of $\mathcal{K} = 15$ nodes (red cross) over the domain $\Omega$ (right).} 
	\label{fig1}
\end{figure}\\
Recalling that the full quadrature is composed by $\mathcal{N}=33152$ points we can say, in the light of the results shown, that our method succeeds to provide a sparse quadrature rule for which the integration error that is introduced can be controlled by specifying $\epsilon$.\\
We again compare our strategy, for $N=7$, with the ones based on $\ell^1$-norm minimization and on the non negative least square method. Table \ref{T5} reports the number of nodes in the sparse rules and the error in the integration of the nonlinear term. As for the other test problem our method is the only one always able to return errors smaller than the required accuracy while the other two approaches fails for the lowest values of $\epsilon$. When $\epsilon<10^{-6}$, the train sample turned out to be not large enough to make the error smaller then the target accuracy. In Figure \ref{fig10} we plot the error with respect to sparsity for a larger train sample of $|\Xi^{\text{rb},\text{train}}|=15^3$ points. With this dataset we indeed manage to recover quadratures of prescribed accuracy $\epsilon$ smaller than $10^{-6}$. Figure \ref{fig9} compares the algorithms in terms of error in the integration of the nonlinear term with respect to sparsity. As in the first test problem the \texttt{FOCUSS} algorithm is the one recovering the most accurate quadrature rule for a fixed number of quadrature points.\\
Last comparison concerns the average computational time required to find the sparse rule for different train samples. Results are displayed in Table \ref{T6} and they reflect what already seen for the Schrödinger fundamental solution test problem. The fastest method is the empirical cubature, while the linear programming is the lowest in any case tested. Also we observe that the computational time related to the empirical cubature strongly increases for the smaller values of $\epsilon$. 
\begin{table}[t]
	\begin{center}
		\begin{tabular}{llllccc}
			\cline{1-7}
			\multirow{3}{*}{$\eps$}&\multirow{3}{*}{$\mathcal{K}_{lp}$}&\multirow{3}{*}{$\mathcal{K}_{ec}$}&\multirow{3}{*}{$\mathcal{K}_{fo}$}&\multirow{3}{*}{$E_{lp}(\Xi^{\text{test}})$}&\multirow{3}{*}{$E_{ec}(\Xi^{\text{test}})$}&\multirow{3}{*}{$E_{fo}(\Xi^{\text{test}})$}\\
			&&&&&&
			\\
			&&&&&&
			\\
			\cline{1-7}
			&&&&&&
			\\
			$10^{-1}$&$6$&8&$40$&$1.0\cdot10^{-1}$&$4.7\cdot10^{-3}$&$1.1\cdot10^{-4}$\\
			$10^{-2}$&$10$&15&$56$&$1.0\cdot10^{-2}$&$1.6\cdot10^{-3}$&$2.9\cdot10^{-5}$\\
			$10^{-3}$&$14$&23&$74$&$1.0\cdot10^{-3}$&$2.8\cdot10^{-4}$&$2.2\cdot10^{-6}$\\
			$10^{-4}$&$21$&36&$92$&$1.0\cdot10^{-3}$&$3.9\cdot10^{-5}$&$1.8\cdot10^{-6}$\\
			$10^{-5}$&$37$&54&$112$&$2.1\cdot10^{-4}$&$1.5\cdot10^{-5}$&$1.4\cdot10^{-6}$\\
			$10^{-6}$&$49$&$67$&$136$&$3.8\cdot10^{-5}$&$9.2\cdot10^{-6}$&$2.4\cdot10^{-7}$\\
		\cline{1-7}
		\end{tabular}   
	\end{center}
	\caption{Parametrized nonlinear diffusion test problem: results in terms of number of quadrature nodes and errors (\ref{eq:err}) for linear programming ($lp$), empirical cubature ($ec$) and \texttt{FOCUSS} based strategy ($fo$).}
	\label{T5}
\end{table}
\begin{figure}[]

	\centering{
		
		\includegraphics[width=0.48\textwidth]{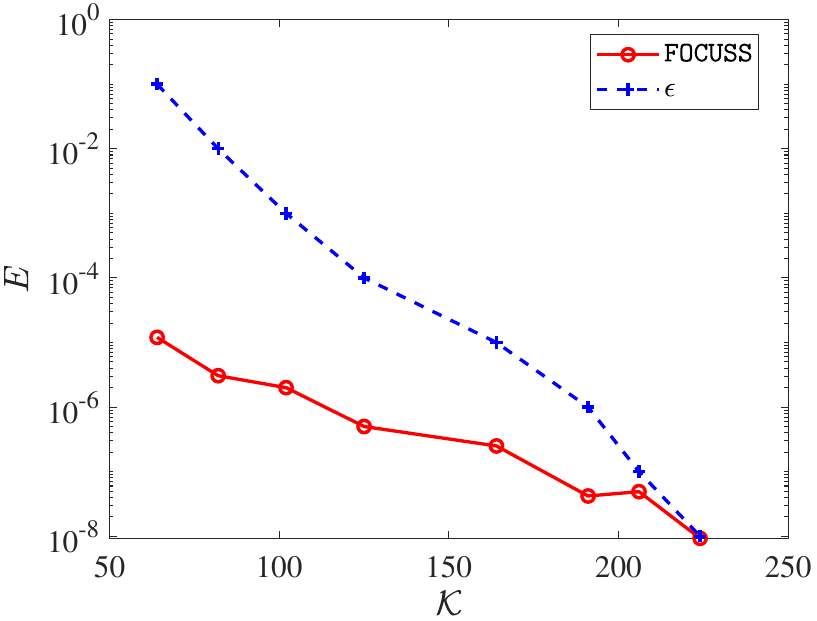}
	}
	
	\caption{Parametrized nonlinear diffusion test problem: error with respect to number of quadrature points for the dataset $|\Xi^{\text{rb},\text{train}}|=15^3$.}
	\label{fig10}
\end{figure}
\begin{figure}[]

	\centering{
		
		\includegraphics[width=0.48\textwidth]{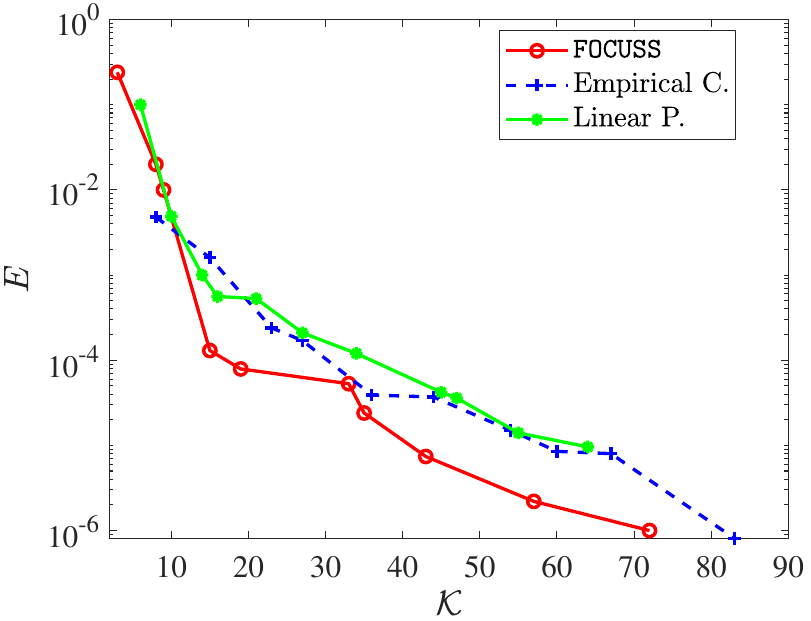}
	}
	
	\caption{Parametrized nonlinear diffusion test problem: error with respect to number of quadrature points for the three sparse recovery algorithms.}
	\label{fig9}
\end{figure}
\begin{table}[t]
	\centering
	\begin{tabular}{|c|ccc|ccc|ccc|}
		\hline
		\multirow{2}{*}{$\epsilon$}&\multicolumn{3}{c|}{$|\Xi^{\text{train}}|=6^3$}&\multicolumn{3}{c|}{$|\Xi^{\text{train}}|=7^3$}&\multicolumn{3}{c|}{$|\Xi^{\text{train}}|=8^3$}\\
		\cline{2-10}
		& $t_{lp}$ &$t_{ec}$& $t_{fo}$ & $t_{lp}$ &$t_{ec}$& $t_{fo}$& $t_{lp}$ &$t_{ec}$& $t_{fo}$\\
		\hline
		\hline
		$10^{-1}$   & 157    &1&17&
		220&2&18&
		308&3&26 \\
		
		$10^{-2}$   & 155    &2&17&
		220&2&18&
		314&5&27 \\
		
		$10^{-3}$    & 155 &2&17&
		240&4&20&
		309&7&30\\
		
		$10^{-4}$    & 157 &3&17&
		226&5&21&
		316&10&31\\
		
		$10^{-5}$  & 156 &4&17&
		218&8&21&
		313&20&30\\
		
		$10^{-6}$  & 156 &8&18&
		223&14&23&
		327&27&30\\
		
		\hline
	\end{tabular}
	\cprotect\caption{Parametrized nonlinear diffusion test problem: average time(s) to compute the sparse rule with linear programming ($lp$), empirical cubature ($ec$) and \texttt{FOCUSS} based method ($fo$) for different sizes of the training dataset $|\Xi^{\text{train}}|$ and different integration tolerance $\epsilon$.}
	\label{T6}
\end{table}
\section{Conclusions}
In this paper, we developed and analyzed an offline/online computational procedure for computing integrals of parametrized functions. The main features are an empirical dataset, from which we extract the relevant information by employing the truncated singular value decomposition, and an empirical quadrature procedure based on an $\ell^p$-quasi-norm minimization problem, which: 1) accommodates the problem's constraints naturally, 2) resolves in a simple numerical scheme and 3) allows efficient calculations in an offline/online setting. We presented theoretical and numerical results to justify our approach. Also, we compared our procedure with the empirical quadrature based on the solution of an $\ell^1$-norm minimization problem and with the one based on the solution of not negative least square problems, both implemented in \texttt{Matlab}. Our method was the only one always able to recover sparse quadrature rule with error smaller or equal to the required tolerance $\epsilon$. Moreover the \texttt{FOUCSS} algorithm showed to provide more accurate rules than the other methods for a fixed number of quadrature points. In terms of time required to compute the empirical rule, the $\ell^p$-norm minimization routine proved to be faster than the $\ell^1$-norm minimization routine. The empirical cubature was the fastest for the higher values of $\epsilon$ but it failed to converge or had poor resolution for the lowest $\epsilon$. Moreover we observed that the computational time demanded by our $\ell^p$-norm minimization routine is dominated by the initial truncated SVD, which cost does not scale linearly with respect to the size of the dataset. Anyway, as already mentioned, for very large datasets one could take advantage of randomized algorithms for matrix decomposition, see \cite{p20}, to strongly speed up the computation related to the truncated SVD.\\
Further developments of this method and specific applications to model order reduction and hyper-reduction will be considered as the object of upcoming publications.

\subsection*{Data availability}

The codes implementing the algorithms discussed in this article are publicly available at:\\
\url{https://github.com/MattiaManucci/Sparse-data-driven-quadrature-rules-via-FOCUSS.git}.

No other data are associated to the manuscript.

\end{document}